\documentclass[preprint,10.5pt]{elsarticle}

\usepackage{amssymb}
\usepackage{amsmath}

\usepackage{mathrsfs}
\usepackage{amsfonts}
\usepackage{CJK,fancybox,color,graphicx,amsmath,amsthm,amssymb,enumerate}

\journal{}

\newtheorem{Lemma}{Lemma}[section]
\newtheorem{Theorem}{Theorem}[section]

\newtheorem{Remark}{Remark}[section]

\theoremstyle{definition}

\theoremstyle{remark}

\numberwithin{equation}{section}

\newcommand{\D}{\displaystyle}
\newcommand{\DF}[2]{\frac{\D#1}{\D#2}}

\begin{document}

\begin{frontmatter}



\title{Concentration in the flux approximation limit of Riemann solutions to the extended Chaplygin gas equations}


\author{Qingling Zhang}

\address{School of Mathematics and Computer Sciences, Jianghan University, Wuhan 430056, PR China }
\ead{zhangqingling2002@163.com}
\begin{abstract}

In this paper, two kinds of occurrence mechanism on the phenomenon
of concentration and the formation of delta shock wave in the flux
approximation limit of Riemann solutions to the extended Chaplygin
gas equations are analyzed and identified. Firstly, the Riemann
problem of the extended Chaplygin gas equations is solved
completely. Secondly, we rigorously show that, as the pressure
vanishes, any two-shock Riemann solution to the extended Chaplygin
gas equations tends to a $\delta$-shock solution to the transport
equations, and the intermediate density between the two shocks tends
to a weighted $\delta$-measure that forms the $\delta$-shock; any
two-rarefaction-wave Riemann solution to the extended Chaplygin gas
equations tends to a two-contact-discontinuity solution to the
transport equations, and the nonvacuum intermediate state between
the two rarefaction waves tends to a vacuum state. At last, we also
show that, as the pressure approaches the generalized Chaplygin
pressure, any two-shock Riemann solution tends to a delta-shock
solution to the generalized Chaplygin gas equations.

\end{abstract}


\end{frontmatter}



\section{Introduction }

\setcounter{equation}{0}

The extended Chaplygin gas equations can be expressed as
\begin{equation}\label{1.1}
\left\{\begin{array}{ll}
\rho_t+(\rho u)_x=0,\\
(\rho u)_t+(\rho u^2+P)_x=0,
\end{array}\right.
\end{equation}
where $\rho$, $u$ and $P$ represent the density, the velocity and
the scalar pressure, respectively, and
\begin{equation}\label{1.2}
P=A\rho^{n}-\frac{B}{\rho^{\alpha}},\ \ \ 1 \leq n \leq 3,\ \
0<\alpha\leq1,
\end{equation}
with two parameters $A,B>0$. This model was proposed by Naji in 2014
\cite{Naji} to study the evolution of dark energy. When $B=0$ in
$(\ref{1.2})$, $P=A\rho^{n}$ is the standard equation of state of
perferct fluid. Up to now, various kinds of theoretical models have
been proposed to interpret the behavior of dark energy. Specially,
when $n=1$ in $(\ref{1.2})$, it reduces to the state equation for
modified Chaplygin gas, which was originally proposed by Benaoum in
2002 \cite{Benaoum}. As an exotic fluid, such a gas can explain the
current accelerated expansion of the universe. Whereas when $A=0$ in
$(\ref{1.2})$, $P=-\frac{B}{\rho^{\alpha}}$ is called the pressure
for the generalized Chaplygin gas \cite{Setare}. Furthermore, when
$\alpha=1$, $P=-\frac{B}{\rho}$ is called the pressure for (pure)
Chaplygin gas which was introduced by Chaplyin \cite{Chaplygin},
Tsien \cite{Tsien} and von Karman \cite{von Karman} as a suitable
mathematical approximation for calculating the lifting force on a
wing of an airplane in aerodynamics. Such a gas own a negative
pressure and occurs in certain theories of cosmology. It has been
also advertised as a possible model for dark energy
\cite{Bilic-Tupper-Viollier,Gorini-Kamenshchik-Moschella-Pasquier,Setare1}.

When two parameters $A$, $B\rightarrow0$, the limit system of
$(\ref{1.1})$ with $(\ref{1.2})$ formally becomes the following
transport equations:
\begin{equation}\label{1.3}
\left\{\begin{array}{ll}
\rho_t+(\rho u)_x=0,\\
(\rho u)_t+(\rho u^2)_x=0,
\end{array}\right.
\end{equation}
which was also called the zero-pressure gas dynamics, and can be
derived from Boltzmann equations \cite{Bouchut} and the
flux-splitting numerical schemes for the full compressible Euler
equations \cite{Brenier-Grenier,Li-Cao}. It can also be used to
describe some important physical phenomena, such as the motion of
free particles sticking together under collision and the formation
of large scale structures in the universe
\cite{Agarwal-Halt,E-Rykov-Sinai,Shandarin-Zeldovich}.

The transport equation $(\ref{1.3})$ has been studied extensively
since 1994. The existence of measure solutions of the Riemann
problem was first proved by Bouchut \cite{Bouchut} and the existence
of the global weak solutions was obtained by Brenier and Grenier
\cite{Bouchut} and E.Rykov and Sinai \cite{E-Rykov-Sinai}. Sheng and
Zhang \cite{Sheng-Zhang} discovered that the $\delta$-shock and
vacuum states do occur in the Riemann solutions to the transport
equation $(\ref{1.3})$ by the vanishing viscosity method. Huang and
Wang \cite{Huang-Wang} proved the uniqueness of the weak solution
when the initial data is a Radon measure. Also see
\cite{Shen1,Wang-Ding,Wang-Huang-Ding,Wang-Zhang1,Yang} for more
related results.

$\delta$-shock is a kind of nonclassical nonlinear waves on which at
least one of the state variables becomes a singular measure.
Korchinski \cite{Korchinski} firstly introduced the concept of the
$\delta$-function into the classical weak solution in his
unpublished Ph. D. thesis. Tan, Zhang and Zheng
\cite{Tan-Zhang-Zheng} considered some 1-D reduced system and
discovered that the form of $\delta$-functions supported on shocks
was used as parts in their Riemann solutions for certain initial
data. LeFloch et al. \cite{LeFloch-Liu} applied the approach of
nonconservative product to consider nonlinear hyperbolic systems in
the nonconservative form. We can also refer to
\cite{Bouchut,Li-Zhang-Yang,Sheng-Zhang} for related equations and
results. Recently, the weak asymptotic method was widely used to
study the $\delta$-shock wave type solution by Danilov and
Shelkovich et
al.\cite{Danilvo-Shelkovich1,Danilvo-Shelkovich2,Shelkovich}.

As for delta shock waves, one research focus is to explore the
phenomena of concentration and cavitation and the formation of delta
shock waves and vacuum states in solutions. In \cite{Chen-Liu1},
Chen and Liu considered the Euler equations for isentropic fluids,
i.e., in $(\ref{1.1})$ they took the prototypical pressure function
as follows:
\begin{equation}\label{1.4}
P=\varepsilon\frac{\rho^\gamma}{\gamma},\ \ \gamma>1.
\end{equation}
They analyzed and identified the phenomena of concentration and
cavitation and the formation of $\delta$-shocks and vacuum states as
$\varepsilon\rightarrow 0$, which checked the numerical observation
for the 2-D case by Chang, Chen and Yang
\cite{Chang-Chen-Yang1,Chang-Chen-Yang2}. They also pointed out that
the occurrence of $\delta$-shocks and vacuum states in the process
of vanishing pressure limit can be regarded as a phenomenon of
resonance between the two characteristic fields. In
\cite{Chen-Liu2}, they made a further step to generalize this result
to the nonisentropic fluids. Specially, for $\gamma=1$ in
$(\ref{1.4})$, the vanishing pressure limit has been studied by Li
\cite{Li}. Besides, the results were extended to the relativistic
Euler equations for polytropic gases by Yin and Sheng
\cite{Yin-Sheng}, the perturbed Aw-Rascle model by Shen and Sun
\cite{Shen-Sun} and the modified Chaplygin gas equations for by Yang
and Wang \cite{Yang-Wang,Yang-Wang1}. For other related works, we
can also see \cite{Mitrovic-Nedeljkov,Zhang}.

In this paper, we focus on the extended Chaplygin gas equations
$(\ref{1.1})$ to discuss the phenomena of concentration and
cavitation and the formation of delta shock waves and vacuum states
in solutions as the double parameter pressure vanishes wholly or
partly, which corresponds to a two parameter limit of solutions in
contrast to the previous works in
\cite{Chen-Liu1,Chen-Liu2,Mitrovic-Nedeljkov,Shen-Sun,Yin-Sheng}.
Equivalently, we will study the limit behavior of Riemann solutions
to the extended Chaplygin gas equations as the pressure vanishes, or
tends to the generalized Chaplygin pressure.

It is noticed that, When $A$, $B\rightarrow0$, the system
$(\ref{1.1})$ with $(\ref{1.2})$ formally becomes the transport
equations $(\ref{1.3})$. For fixed $B$, When $A\rightarrow0$, the
system $(\ref{1.1})$ with $(\ref{1.2})$ formally becomes the
generalized Chaplygin gas equations
\begin{equation}\label{1.5}
\left\{\begin{array}{ll}
\rho_t+(\rho u)_x=0,\\
(\rho u)_t+(\rho u^2-\frac{B}{\rho^{\alpha}})_x=0.
\end{array}\right.
\end{equation}
 When $\alpha=1$, it is just the Chaplygin gas equations.
In 1998, Brenier \cite{Brenier} firstly studied the 1-D Riemann
problem and obtained the solutions with concentration when initial
data belong to a certain domain in the phase plane. Recently, Guo,
Sheng and Zhang \cite{Guo-Sheng-Zhang} abandoned this constrain and
constructively obtained the general solutions of the 1-D Riemann
problem in which the $\delta$-shock wave developed. Moreover, in
that paper, they also systematically studied the 2-D Riemann problem
for isentropic Chaplygin gas equations. In \cite{Wang}, Wang solved
the Riemann problem of $(\ref{1.5})$ by the weak asymptotic method.
It has been shown that, in their results, $\delta$-shocks do occur
in the Riemann solutions, but vacuum states do not. For more results
about Chaplygin gas, one can refer to
\cite{Qu-Wang,Shen2,Sun,Wang-Zhang2}.

In this paper, we first solve the Riemann problem of system
(\ref{1.1}) with Riemann initial data
\begin{equation}\label{1.6}
(\rho,u)(x,0)=(\rho_\pm,u_\pm),\ \ \ \ \pm x>0,
\end{equation}
where $\rho_\pm>0,\ u_\pm$ are arbitrary constants. With the help of
analysis method in phase plane, we constructed the Riemann solutions
with four different structures: $R_1R_2$, $R_1S_2$, $S_1R_2$ and
$S_1S_2$.

Then we analyze the formation of $\delta$-shocks and vacuum states
in the Riemann solutions as the pressure vanishes. It is shown that,
as the pressure vanishes, any two-shock Riemann solution to the
extended Chaplygin gas equations tends to a $\delta$-shock solution
to the transport equations, and the intermediate density between the
two shocks tends to a weighted $\delta$-measure that forms the
$\delta$-shock; by contrast, any two-rarefaction-wave Riemann
solution to the extended Chaplygin gas equations tends to a
two-contact-discontinuity solution to the transport equations, and
the nonvacuum intermediate state between the two rarefaction waves
tends to a vacuum state, even when the initial data stay away from
the vacuum. As a result, the delta shocks for the transport
equations result from a phenomenon of concentration, while the
vacuum states results from a phenomenon of cavitation in the
vanishing pressure limit process. These results are completely
consistent with that in \cite{Chen-Liu1}, and also cover those
obtained in \cite{Yang-Wang,Yang-Wang1}.

In addition, we also proved that as the pressure tends to the
generalized Chaplygin pressure ($A\rightarrow0$), any two-shock
Riemann solution to the extended Chaplygin gas equations tends to a
$\delta$-shock solution to the generalized Chaplygin gas equations,
and the intermediate density between the two shocks tends to a
weighted $\delta$-measure that forms the $\delta$-shock; by
contrast, any two-rarefaction-wave Riemann solution to the extended
Chaplygin gas equations tends to the two-rarefaction-wave
(two-contact-discontinuity for $\alpha=1$) solution to the transport
equations, and the intermediate state between the two rarefaction
waves (two contact discontinuities) is a nonvacuum state.
Consequently, the delta shocks for the  generalized Chaplygin gas
equations result from a phenomenon of concentration in the partly
vanishing pressure limit process.

From the above analysis, we can find two kinds of occurrence
mechanism on the phenomenon of concentration and the formation of
delta shock wave. On one hand, since the strict hyperbolicity of the
limiting system $(\ref{1.3})$ fails, see Section 4, the delta shock
wave forms in the limit process as the pressure vanishes. This is
consistent with those results obtained in
\cite{Chen-Liu1,Chen-Liu2,Mitrovic-Nedeljkov,Shen-Sun,Yang-Wang,Yin-Sheng}.
On the other hand, the strict hyperbolicity of the limiting system
$(\ref{1.5})$ is preserved, see Section 5, the formation of delta
shock waves still occur as the pressure partly vanishes. In this
regard, it is different from that in
\cite{Chen-Liu1,Chen-Liu2,Mitrovic-Nedeljkov,Shen-Sun,Yang-Wang,Yin-Sheng}.
In any case, the phenomenon of concentration and the formation of
delta shock wave can be regarded as a process of resonance formation
between two characteristic fields.

The paper is organized as follows. In Section 2, we restate the
Riemann solutions to transport equations $(\ref{1.3})$ and the
generalized Chaplygin gas equations $(\ref{1.5})$. In Section 3, we
investigate the Riemann problem of the extended Chaplygin gas
equations $(\ref{1.1})$-$(\ref{1.2})$ and examine the dependence of
the Riemann solutions on the two parameters $A,B>0$. In Section 4,
we analyze the limit of Riemann solutions to the extended Chaplygin
gas equations $(\ref{1.1})$-$(\ref{1.2})$ with $(\ref{1.6})$ as the
pressure vanishes. In Section 5, we discuss the limit of Riemann
solutions to the extended Chaplygin gas equations
$(\ref{1.1})$-$(\ref{1.2})$ with $(\ref{1.6})$ as the pressure
approaches to the generalized Chaplygin pressure. Finally,
conclusions and discussions are drawn in Section 6.

\section{Preliminaries}
\subsection{Riemann problem for the transport equations} In this
section, we restate the Riemann solutions to the transport equations
$(\ref{1.3})$ with initial data $(\ref{1.6})$. See
\cite{Sheng-Zhang} for more details.

The transport equations $(\ref{1.3})$ have a double eigenvalue
$\lambda=u$ and only one right eigenvectors $ \vec{r}=(1,0)^T $.
Furthermore, we have $\nabla\lambda\cdot\vec{r}=0$, which means that
$\lambda$ is linearly degenerate. The Riemann problem $(\ref{1.3})$
and $(\ref{1.6})$ can be solved by contact discontinuities, vacuum
or $\delta$-shocks connecting two constant states
$(\rho_\pm,u_\pm)$.

By taking the self-similar transformation $\xi=\frac{x}{t}$, the
Riemann problem is reduced to the boundary value problem of the
ordinary differential equations:
\begin{equation}\label{2.1}
\left\{\begin{array}{ll}
-\xi\rho_\xi+(\rho u)_\xi=0,\\
-\xi(\rho u)_\xi+(\rho u^2)_\xi=0,
\end{array}\right.
\end{equation}
with $(\rho,u)(\pm\infty)=(\rho_\pm,u_\pm)$.

For the case $u_-<u_+$, there is no characteristic passing through
the region $\{\xi: u_-<\xi<u_+\}$, so the vacuum should appear in
the region. The solution can be expressed as
\begin{equation}\label{2.2}
(\rho,u)(\xi)=\left\{\begin{array}{ll}
(\rho_-,u_-),\ \ \ \ -\infty<\xi\leq u_-,\\
(0,\xi),\ \ \ \ \ \ \ \ u_-<\xi< u_+,\\
(\rho_+,u_+),\ \ \ \ u_+\leq\xi<\infty.
\end{array}\right.
\end{equation}

For the case $u_-=u_+$, it is easy to see that the constant states
$(\rho_\pm,u_\pm)$ can be connected by a contact discontinuity.

For the case $u_->u_+$, a solution containing a weighted
$\delta$-measure supported on a curve will be constructed. Let
$x=x(t)$ be a discontinuity curve, we consider a piecewise smooth
solution of $(\ref{1.3})$ in the form
\begin{equation}\label{2.3}
(\rho,u)(x,t)=\left\{\begin{array}{ll}
(\rho_-,u_-),\ \ \ \ \ \ \ \ \ \ \ \ \ \ \ \ \ \ \ \ \ \ \ \ x<x(t),\\
(w(t)\delta(x-x(t)),u_\delta(t)),\ \ \ \ \ x=x(t),\\
(\rho_+,u_+),\ \ \ \ \ \ \ \ \ \ \ \ \ \ \ \ \ \ \ \ \ \ \ \ x>x(t).
\end{array}\right.
\end{equation}

In order to define the measure solution as above, like as in
\cite{Chen-Liu1,Chen-Liu2,Sheng-Zhang}, the two-dimensional weighted
$\delta$-measure $w(t)\delta_S$ supported on a smooth curve
$S=\{(x(s),t(s)):a\leq s\leq b\}$ should be introduced as follows:
\begin{equation}\label{2.4}
\langle
w(\cdot)\delta_S,\psi(\cdot,\cdot)\rangle=\int_a^bw(s)\psi(x(s),t(s))\sqrt{{x'(s)}^2+{t'(s)}^2}ds,
\end{equation}
for any $\psi\in C_0^\infty(R\times R_{+})$.

As shown in \cite{Sheng-Zhang}, for any $\psi\in C_0^\infty(R\times
R_{+})$, the $\delta$-measure solution $(\ref{2.3})$ constructed
above satisfies
\begin{equation}\label{2.5}
\left\{\begin{array}{ll} \langle\rho,\psi_t\rangle+\langle\rho
u,\psi_x\rangle=0,
\\
\langle\rho u,\psi_t\rangle+\langle\rho u^2,\psi_x\rangle=0,
\end{array}\right.
\end{equation}
 in which
$$\langle\rho,\psi\rangle=\int_0^\infty\int_{-\infty}^\infty\rho_0\psi dxdt+\langle w_1(\cdot)\delta_S,\psi(\cdot,\cdot)\rangle,$$
$$\langle\rho u,\psi\rangle=\int_0^\infty\int_{-\infty}^\infty\rho_0 u_0\psi dxdt+\langle w_2(\cdot)\delta_S,\psi(\cdot,\cdot)\rangle,$$
where
$$\rho_0=\rho_-+[\rho]H(x-\sigma t),\ \ \rho_0u_{0}=\rho_-u_{-}+[\rho u]H(x-\sigma t),$$
and
$$w_1(t)=\frac{t}{\sqrt{1+\sigma^2}}(\sigma[\rho]-[\rho u]), \ \ w_2(t)=\frac{t}{\sqrt{1+\sigma^2}}(\sigma[\rho u]-[\rho u^2]).$$
Here $H(x)$ is the Heaviside function given by $H(x)=1$ for $x>0$
and $H(x)=0$ for $x<0$.

Substituting $(\ref{2.3})$ into $(\ref{2.5})$, one can derive the
generalized Rankine-Hugoniot conditions
\begin{equation}\label{2.6}
\left\{\begin{array}{ll}
\DF{dx(t)}{dt}=u_\delta(t),\\
\DF{dw(t)}{dt}=[\rho] u_\delta(t)-[\rho u],\\
\DF{d(w(t) u_ \delta(t))}{dt}=[\rho u] u_\delta(t)-[\rho u^2]
\end{array}\right.
\end{equation}
where $[\rho]=\rho_+-\rho_-$, etc.

Through solving $(\ref{2.6})$ with $x(0)=0,\ w(t)=0$, we obtain
\begin{equation}\label{2.7}
\left\{\begin{array}{ll}
u_\delta(t)=\sigma=\DF{\sqrt{\rho_-}u_-+\sqrt{\rho_+}u_+}{\sqrt{\rho_-}+\sqrt{\rho_+}},\\
x(t)=\sigma t,\\
w(t)=-\sqrt{\rho_-\rho_+}(u_+-u_-)t,
\end{array}\right.
\end{equation}

Moreover, the $\delta$-measure solution $(\ref{2.3})$ with
$(\ref{2.6})$ satisfies the $\delta$-entropy condition:
\begin{equation}\nonumber
u_+<\sigma<u_-,
\end{equation}
which means that all the characteristics on both sides of the
$\delta$-shock are incoming.

\subsection{Riemann problem for the generalized Chaplygin gas
equations}

In this section, we solve the Riemann problem for the generalized
Chaplygin gas equations $(\ref{1.5})$ with $(\ref{1.6})$,  which one
can also see in \cite{Guo-Sheng-Zhang,Wang}.

It is easy to see that $(\ref{1.5})$ has two eigenvalues
$$\lambda_1^B=u-\sqrt{\alpha B}\rho^{-\frac{\alpha+1}{2}},\ \ \lambda_2^B=u+\sqrt{\alpha B}\rho^{-\frac{\alpha+1}{2}},$$
with corresponding right eigenvectors
$$\overrightarrow{r_1}^B=(-\sqrt{\alpha B}\rho^{-\frac{\alpha+1}{2}},\rho)^T,\ \ \overrightarrow{r_2}^B=(\sqrt{\alpha B}\rho^{-\frac{\alpha+1}{2}},\rho)^T.$$
So $(\ref{1.5})$ is strictly hyperbolic for $\rho>0$. Moreover, when
$0<\alpha<1$, we have
 $\bigtriangledown\lambda_i^B\cdot\overrightarrow{r_i}^B\neq0$,
$i=1,2$, which implies that $\lambda_1^B$ and $\lambda_2^B$ are both
genuinely nonlinear and the associated waves are rarefaction waves
and shock waves. When $\alpha=1$, $\bigtriangledown\lambda_i^B\cdot
\overrightarrow{r_i}^B=0$, $i=1,2$, which implies that $\lambda_1^B$
and $\lambda_2^B$ are both linearly degenerate and the associated
waves are both contact discontinuities, see \cite{Smoller}.

Since system $(\ref{1.5})$ and the Riemann initial data
$(\ref{1.6})$ are invariant under stretching of coordinates
$(x,t)\rightarrow(\beta x,\beta t)$ ($\beta$ is constant), we seek
the self-similar solution
$$(\rho,u)(x,t)=(\rho,u)(\xi),\ \ \xi=\frac{x}{t}.$$
Then the Riemann problem $(\ref{1.5})$ and $(\ref{1.6})$ is reduced
to the following boundary value problem of the ordinary differential
equations:
\begin{equation}\label{2.8}
\left\{\begin{array}{ll}
-\xi\rho_\xi+(\rho u)_\xi=0,\\
-\xi(\rho u)_\xi+(\rho u^2-\frac{B}{\rho^\alpha})_\xi=0,\\
\end{array}\right.
\end{equation}
with $(\rho,u)(\pm\infty)=(\rho_\pm,u_\pm)$.

Besides the constant solution, it provides the backward rarefaction
wave
\begin{equation}\label{2.9}
 \overleftarrow{R}(\rho_-,u_-):\left\{\begin{array}{ll} \xi=\lambda_1^B=u-\sqrt{\alpha B}\rho^{-\frac{\alpha+1}{2}},\\
u-\frac{2\sqrt{\alpha
B}}{1+\alpha}\rho^{-\frac{\alpha+1}{2}}=u_--\frac{2\sqrt{\alpha
B}}{1+\alpha}\rho_-^{-\frac{\alpha+1}{2}},\ \ \rho<\rho_-,
\end{array}\right.
\end{equation}
and the forward rarefaction wave
\begin{equation}\label{2.10}
\overrightarrow{R}(\rho_-,u_-):\left\{\begin{array}{ll} \xi=\lambda_2^B=u+\sqrt{\alpha B}\rho^{-\frac{\alpha+1}{2}},\\
u+\frac{2\sqrt{\alpha
B}}{1+\alpha}\rho^{-\frac{\alpha+1}{2}}=u_-+\frac{2\sqrt{\alpha
B}}{1+\alpha}\rho_-^{-\frac{\alpha+1}{2}},\ \ \rho>\rho_-,
\end{array}\right.
\end{equation}
When $\alpha=1$, the backward (forward) rarefaction wave becomes the
backward (forward) contact discontinuity.

For a bounded discontinuity at $\xi=\sigma$, the Rankine-Hugoniot
conditions hold:
\begin{equation}\label{2.11}
\left\{\begin{array}{ll}
-\sigma^{B}[\rho]+[\rho u]=0,\\
-\sigma^{B}[\rho u]+[\rho u^2-\DF{B}{\rho^\alpha}]=0,
\end{array}\right.
\end{equation}
where $[\rho]=\rho-\rho_-$, etc. Together with the Lax shock
inequalities, (\ref{2.11}) gives the backward shock wave
\begin{equation}\label{2.12}
\overleftarrow{S}(\rho_-,u_-):\left\{\begin{array}{ll}
\sigma_1^{B}=\DF{\rho u-\rho_-u_-}{\rho-\rho_-},\\
u-u_-=-\sqrt{B(\frac{1}{\rho}-\frac{1}{\rho_-})(\frac{1}{\rho^{\alpha}}-\frac{1}{\rho_{-}^{\alpha}})},\
\ \rho>\rho_-,
\end{array}\right.
\end{equation}
and the forward shock wave
\begin{equation}\label{2.13}
\overrightarrow{S}(\rho_-,u_-):\left\{\begin{array}{ll}
\sigma_2^{B}=\DF{\rho u-\rho_-u_-}{\rho-\rho_-},\\
u-u_-=-\sqrt{B(\frac{1}{\rho}-\frac{1}{\rho_-})(\frac{1}{\rho^{\alpha}}-\frac{1}{\rho_{-}^{\alpha}})},\
\ \rho<\rho_-.
\end{array}\right.
\end{equation}
When $\alpha=1$, the backward (forward) shock wave becomes the
backward (forward) contact discontinuity.

Furthermore, for a given left state $(\rho_-,u_-)$, the backward
shock wave $\overleftarrow{S}(\rho_-,u_-)$ has a straight line
$u=u_--\sqrt{ B}\rho_-^{-\frac{\alpha+1}{2}}$, as its asymptote, and
for a given right state $(\rho_+,u_+)$, the forward shock wave
$\overrightarrow{S}(\rho_+,u_+)$ has a straight line $u=u_++\sqrt{
B}\rho_+^{-\frac{\alpha+1}{2}}$ as its asymptote.

It is easy to see that, when $u_++\sqrt{
B}\rho_+^{-\frac{\alpha+1}{2}} \leq u_--\sqrt{
B}\rho_-^{-\frac{\alpha+1}{2}}$, the backward shock wave
$\overleftarrow{S}(\rho_-,u_-)$ can not intersect the forward shock
wave $\overrightarrow{S}(\rho_+,u_+)$, a delta shock wave must
develop in solutions. Under the definition (\ref{2.4}), a delta
shock wave can be introduced to construct the solution of
(\ref{1.5})-(\ref{1.6}), which can be expressed as
\begin{equation}\label{2.14}
(\rho,u)(x,t)=\left\{\begin{array}{ll}
(\rho_-,u_-),\ \ \ \ \ \ \ \ \ \ \ \ \ \ \ \ \ \ \ \ \ \ \ x<\sigma^B t,\\
(w^B(t)\delta(x-\sigma^B t),\sigma^B),\ \ \  \ \ x=\sigma^B t,\\
(\rho_+,u_+),\ \ \ \ \ \ \ \ \ \ \ \ \ \ \ \ \ \ \ \ \ \ \
x>\sigma^B t,
\end{array}\right.
\end{equation}
with
\begin{equation}\nonumber
\frac{B}{\rho^\alpha}=\left\{\begin{array}{ll}
\DF{B}{\rho_-^\alpha},\ \ \ \ \ \ \ x<\sigma^B t,\\
0,\ \ \  \ \ \ \ \ \ x=\sigma^B t,\\
\DF{B}{\rho_+^\alpha},\ \ \ \ \ \ \ x>\sigma^B t,
\end{array}\right.
\end{equation} see
\cite{Brenier}.

By the weak solution definition in Subsection 2.1, for the system
(\ref{1.5})£¬we can get the following generalized Rankine-Hugoniot
conditions
\begin{equation}\label{2.15}
\left\{\begin{array}{ll}
\DF{dx^{B}(t)}{dt}=u_\delta^{B}(t)=\sigma^B,\\[4pt]
\DF{dw^{B}(t)}{dt}=u_\delta^{B}(t)[\rho]-[\rho u],\\[4pt]
\DF{d(w^{B}(t)u_\delta^{B}(t))}{dt}=u_\delta^{B}(t)[\rho u]-[\rho
u^2-\DF{1}{\rho}],
\end{array}\right.
\end{equation}
where $x^{B}(t)$, $w^{B}(t)$ and $u_\delta^{B}(t)$ are respectively
denote the location, weight and propagation speed of the
$\delta$-shock, $[\rho]=\rho(x^{B}(t)+0,t)-\rho(x^{B}(t)-0,t)$
denotes the jump of the function $\rho$ across the $\delta$-shock.

Then by solving (\ref{2.15}) with initial data $x(0)=0,\ w^B(0)=0$,
under the entropy condition
\begin{equation}\label{2.16}
u_++\sqrt{\alpha B}\rho_+^{-\frac{\alpha+1}{2}}<\sigma^B<
u_--\sqrt{\alpha  B}\rho_-^{-\frac{\alpha+1}{2}},
\end{equation}
 we can obtain
\begin{equation}\label{2.17}
w^B(t)=\big\{\rho_+\rho_-\big((u_+-u_-)^2-(\frac{1}{\rho_+}-\frac{1}{\rho_-})(\frac{B}{\rho_+^\alpha}-\frac{B}{\rho_-^\alpha})\big)\big\}^\frac{1}{2}t,
\end{equation}
\begin{equation}\label{2.18}
\sigma^B=\frac{\rho_+
u_+-\rho_-u_-+\frac{dw^{B}(t)}{dt}}{\rho_+-\rho_-},
\end{equation}
when $\rho_+\neq\rho_-$, and
\begin{equation}\label{2.19}
w^B(t)=(\rho_-u_--\rho_+ u_+)t,
\end{equation}
\begin{equation}\label{2.20}
\sigma^B=\frac{1}{2}(u_++u_-),
\end{equation}
 when $\rho_+=\rho_-$.

In the phase plane ($\rho>0$, $u\in R$), given a constant state
$(\rho_-,u_-)$, we draw the elementary wave curves
(\ref{2.9})-(\ref{2.10}) and (\ref{2.12})-(\ref{2.13}) passing
through this point, which are denoted by
$\overleftarrow{R}$,$\overrightarrow{R}$,$\overleftarrow{S}$ and
$\overrightarrow{S}$ respectively. The backward shock wave
$\overleftarrow{S}$ has an asymptotic line  $u=u_--\sqrt{
B}\rho_-^{-\frac{\alpha+1}{2}}$. In addition, we draw a $S_{\delta}$
curve, which is determined by
\begin{equation}\label{2.21}
u+\sqrt{B}\rho^{-\frac{\alpha+1}{2}}= u_--\sqrt{
B}\rho_-^{-\frac{\alpha+1}{2}}, \ \rho>0.
\end{equation}
Then the phase plane can be divided into five parts I$(\rho_-,u_-)$,
I\!I$(\rho_-,u_-)$, I\!I\!I$(\rho_-,u_-)$, I\!V$(\rho_-,u_-)$ and
V$(\rho_-,u_-)$, see Fig.1.

By the analysis method in the phase plane, one can construct the
Riemann solutions for any given $(\rho_+,u_+)$ as follows:\\
(1) $(\rho_+,u_+)\in$\rm{ I}$(\rho_-,u_-)$:\
$\overleftarrow{R}+\overrightarrow{R}$;\  \ \ \ \ (2)
$(\rho_+,u_+)\in$ I\!I$(\rho_-,u_-)$:\
$\overleftarrow{R}+\overrightarrow{S}$;\\(3) $(\rho_+,u_+)\in$
I\!I\!I$(\rho_-,u_-)$:\ $\overleftarrow{S}+\overrightarrow{R}$;\  \
\ \  (4) $(\rho_+,u_+)\in$ I\!V$(\rho_-,u_-)$:\
$\overleftarrow{S}+\overrightarrow{S}$;\\(5) $(\rho_+,u_+)\in$
V$(\rho_-,u_-)$:\ $\delta$-shock.

\unitlength 1mm 
\linethickness{0.4pt}
\ifx\plotpoint\undefined\newsavebox{\plotpoint}\fi 
\begin{picture}(115,72)(-8,0)
\put(110.75,12){\vector(1,0){.07}} \put(13,12){\line(1,0){97.75}}
\put(9.25,72){\vector(0,1){.07}} \put(9.25,15){\line(0,1){57}}

\bezier{50}(52,68)(52,50)(52,12) \bezier{50}(90,68)(90,50)(90,12)

\qbezier(89.5,67.5)(78.13,15)(39.25,13.5)
\qbezier(52.5,67.5)(60.13,13.5)(103.25,14.5)
\qbezier(50.75,67)(43,14.5)(14.25,13) \tiny{
\put(113,11.5){\makebox(0,0)[cc]{$u$}}
\put(7,70.25){\makebox(0,0)[cc]{$\rho$}}
\put(51.5,8){\makebox(0,0)[cc]{$u_--\sqrt{
B}\rho_-^{-\frac{\alpha+1}{2}}$}}
\put(89.25,8){\makebox(0,0)[cc]{$u_-+\sqrt{
B}\rho_-^{-\frac{\alpha+1}{2}}$}}
\put(94.75,33.25){\makebox(0,0)[cc]{I$(\rho_-,u_-)$}}
\put(71.75,16.75){\makebox(0,0)[cc]{I\!I$(\rho_-,u_-)$}}
\put(70.25,44.75){\makebox(0,0)[cc]{I\!I\!I$(\rho_-,u_-)$}}
\put(70,26.5){\circle*{.6}} \put(37.5,27){\circle*{.6}}
\put(23,27){\makebox(0,0)[cc]{$(\rho_-,u_--2\sqrt{
B}\rho_-^{-\frac{\alpha+1}{2}})$}}
\put(81,26.5){\makebox(0,0)[cc]{$(\rho_-,u_-)$}}
\put(48.25,30.25){\makebox(0,0)[cc]{I\!V$(\rho_-,u_-)$}}
\put(25.75,37.5){\makebox(0,0)[cc]{V$(\rho_-,u_-)$}}
\put(58,58.5){\makebox(0,0)[cc]{$\overleftarrow{S}$}}
\put(84,58.25){\makebox(0,0)[cc]{$\overrightarrow{R}$}}
\put(95,17.25){\makebox(0,0)[cc]{$\overleftarrow{R}$}}
\put(50,17){\makebox(0,0)[cc]{$\overrightarrow{S}$}}
\put(42,44.5){\makebox(0,0)[cc]{$S_{\delta}$}}
\put(52,2){\makebox(0,0)[cc]{Fig.1}}}
\end{picture}


\section{Riemann problem for the extended Chaplygin gas equations }
In this section, we first solve the elementary waves and construct
solutions to the Riemann problem of $(\ref{1.1})$-$(\ref{1.2})$ with
$(\ref{1.6})$, and then examine the dependence of the Riemann
solutions on the two parameters $A,B>0$.

The eigenvalues of the system $(\ref{1.1})$-$(\ref{1.2})$ are
$$\lambda_1^{AB}=u-\sqrt{An\rho^{n-1}+\frac{\alpha B}{\rho^{\alpha+1}}}, \quad \lambda_2^{AB}=u+\sqrt{An\rho^{n-1}+\frac{\alpha B}{\rho^{\alpha+1}}},$$
with corresponding right eigenvectors
$$\vec{r}_1^{AB}=(-\rho,\sqrt{An\rho^{n-1}+\frac{\alpha B}{\rho^{\alpha+1}}})^T,\quad \vec{r}_2^{AB}=(\rho,\sqrt{An\rho^{n-1}+\frac{\alpha B}{\rho^{\alpha+1}}})^T.$$

Moreover, we have
$$\nabla\lambda_i^{AB}\cdot \vec{r}_i^{AB}=\DF{An(n+1)\rho^{n+\alpha}+(1-\alpha)\alpha B}{2\sqrt{(An\rho^{n+\alpha}+\alpha B)\rho^{\alpha+1}}}>0\ \ (i=1,2).$$
Thus $\lambda_1^{AB}$ and $\lambda_2^{AB}$ are genuinely nonlinear
and the associated elementary waves are shock waves and rarefaction
waves.

For $(\ref{1.1})$-$(\ref{1.2})$ with $(\ref{1.6})$ are invariant
under uniform stretching of coordinates:$(x,t)\rightarrow(\beta
x,\beta t)$ where constant $\beta>0$, we seek the self-similar
solution$$(\rho,u)(x,t)=(\rho(\xi),u(\xi)),\ \ \xi=\frac{x}{t}.$$
Then the Riemann problem $(\ref{1.1})$-$(\ref{1.2})$ with
$(\ref{1.6})$ is reduced to the boundary value problem of the
following ordinary differential equations:
\begin{equation}\label{3.1}
\left\{\begin{array}{ll}
-\xi\rho_\xi+(\rho u)_\xi=0,\\
-\xi(\rho u)_\xi+(\rho u^2+P)_\xi=0,\ \
P=A\rho^{n}-\frac{B}{\rho^{\alpha}},
\end{array}\right.
\end{equation}
with $(\rho,u)(\pm\infty)=(\rho_\pm,u_\pm)$.

Any smooth solutions of $(\ref{3.1})$ satisfies
\begin{equation}\label{3.2}
\left(\begin{array}{lll}
u-\xi &\ \  \rho  \\
An\rho^{n-1}+\frac{\alpha B}{\rho^{\alpha+1}} &\ \  \rho(u-\xi)
\end{array}\right)
\left(\begin{array}{lll}
 d\rho\\du
\end{array}\right)=0.
\end{equation}
It provides either the constant state solutions
$$(\rho,u)(\xi)=\rm{constant}, $$
or the rarefaction wave which is a continuous solutions of
$(\ref{3.2})$ in the form $(\rho,u)(\xi)$. Then, according to
\cite{Smoller}, for a given left state $(\rho_-,u_-)$, the
rarefaction wave curves in the phase plane, which are the sets of
states that can be connected on the right by a 1-rarefaction wave or
2-rarefaction wave, are as follows:
\begin{equation}\label{3.3}
R_1(\rho_-,u_-): \left\{\begin{array}{ll} \xi=\lambda_1=u-\sqrt{An\rho^{n-1}+\frac{\alpha B}{\rho^{\alpha+1}}},\\
u-u_-=-\int_{\rho_{-}}^{\rho}\frac{\sqrt{An\rho^{n-1}+\frac{\alpha
B}{\rho^{\alpha+1}}}}{\rho}d\rho,
\end{array}\right.
\end{equation}
and
\begin{equation}\label{3.4}
R_2(\rho_-,u_-): \left\{\begin{array}{ll} \xi=\lambda_2=u+\sqrt{An\rho^{n-1}+\frac{\alpha B}{\rho^{\alpha+1}}},\\
u-u_-=\int_{\rho_{-}}^{\rho}\frac{\sqrt{An\rho^{n-1}+\frac{\alpha
B}{\rho^{\alpha+1}}}}{\rho}d\rho.
\end{array}\right.
\end{equation}

From $(\ref{3.3})$ and $(\ref{3.4})$, we obtain that

\begin{eqnarray}
\frac{d\lambda_1^{AB}}{d\rho}=\frac{\partial\lambda_1^{AB}}{\partial
u}\frac{du}{d\rho}+\frac{\partial\lambda_1^{AB}}{\partial\rho}
=-\DF{An(n+1)\rho^{n-1}+\frac{\alpha(1-\alpha)
B}{\rho^{\alpha+1}}}{2\rho\sqrt{An\rho^{n-1}+\frac{\alpha
B}{\rho^{\alpha+1}}}}<0\label{3.5}
\end{eqnarray}

\begin{eqnarray}
\frac{d\lambda_2^{AB}}{d\rho}=\frac{\partial\lambda_2^{AB}}{\partial
u}\frac{du}{d\rho}+\frac{\partial\lambda_2^{AB}}{\partial\rho}
=\DF{An(n+1)\rho^{n-1}+\frac{\alpha(1-\alpha)
B}{\rho^{\alpha+1}}}{2\rho\sqrt{An\rho^{n-1}+\frac{\alpha
B}{\rho^{\alpha+1}}}}>0\label{3.6}
\end{eqnarray}
which imply that the velocity of 1-rarefaction (2-rarefaction) wave
$\lambda_1^{AB}$ ($\lambda_2^{AB}$) is monotonic decreasing
(increasing) with respect to $\rho$.

With the requirement
$\lambda_1^{AB}(\rho_{-},u_{-})<\lambda_1^{AB}(\rho,u)$ and
$\lambda_2^{AB}(\rho_{-},u_{-})<\lambda_2^{AB}(\rho,u)$, noticing
$(\ref{3.5})$ and $(\ref{3.6})$, we get that

\begin{equation}\label{3.7}
R_1(\rho_-,u_-): \left\{\begin{array}{ll} \xi=\lambda_1=u-\sqrt{An\rho^{n-1}+\frac{\alpha B}{\rho^{\alpha+1}}},\\
u-u_-=-\int_{\rho_{-}}^{\rho}\frac{\sqrt{An\rho^{n-1}+\frac{\alpha
B}{\rho^{\alpha+1}}}}{\rho}d\rho,\ \  \rho<\rho_{-},
\end{array}\right.
\end{equation}
and
\begin{equation}\label{3.8}
R_2(\rho_-,u_-): \left\{\begin{array}{ll} \xi=\lambda_2=u+\sqrt{An\rho^{n-1}+\frac{\alpha B}{\rho^{\alpha+1}}},\\
u-u_-=\int_{\rho_{-}}^{\rho}\frac{\sqrt{An\rho^{n-1}+\frac{\alpha
B}{\rho^{\alpha+1}}}}{\rho}d\rho,\ \  \rho>\rho_{-}.
\end{array}\right.
\end{equation}

For the 1-rarefaction wave, through differentiating $u$ respect to
$\rho$ in the second equation in $(\ref{3.7})$, we get
\begin{equation}\label{3.9}
u_{\rho}=-\DF{\sqrt{An\rho^{n-1}+\frac{\alpha
B}{\rho^{\alpha+1}}}}{\rho}<0.
\end{equation}

\begin{equation}\label{3.10}
u_{\rho\rho}=\DF{-An(n-3)\rho^{n+\alpha}+\alpha(\alpha+3)
B}{2\rho^{2}\sqrt{An\rho^{n+\alpha}+\alpha B\rho^{\alpha+1}}}.
\end{equation}

Thus, it is easy to get $u_{\rho\rho}>0$ for $1\leq n\leq3$, i.e.,
the 1-rarefaction wave is convex for $1\leq n\leq3$ in the upper
half phase plane ($\rho>0$).

In addition, from the second equation of $(\ref{3.7})$, we have
$$u-u_-=\int_{\rho}^{\rho_{-}}\frac{\sqrt{An\rho^{n-1}+\frac{\alpha
B}{\rho^{\alpha+1}}}}{\rho}d\rho\geq\int_{\rho}^{\rho_{-}}\sqrt{\alpha
B}\rho^{-\frac{\alpha+1}{2}-1}d\rho=\frac{2\sqrt{\alpha
B}}{\alpha+1}(\rho^{-\frac{\alpha+1}{2}}-\rho_{-}^{-\frac{\alpha+1}{2}}),$$
which means that $\lim\limits_{\rho\rightarrow0}u=+\infty.$

By a similar computation, we have that, for the 2-rarefaction wave,
$u_{\rho}>0$, $u_{\rho\rho}<0$ for $1\leq n\leq3$ and
$\lim\limits_{\rho\rightarrow+\infty}u=+\infty.$ Thus, we can draw
the conclusion that the 2-rarefaction wave is concave for $1\leq
n\leq3$ in the upper half phase plane ($\rho>0$).

Now we consider the discontinuous solution. For a bounded
discontinuity at $\xi=\sigma$, the Rankine-Hugoniot condition holds:
\begin{equation}\label{3.11}
\left\{\begin{array}{ll}
\sigma^{AB}[\rho]=[\rho u],\\
\sigma^{AB}[\rho u]=[\rho u^2+P],\ \
P=A\rho^{n}-\frac{B}{\rho^{\alpha}},
\end{array}\right.
\end{equation}
where $[\rho]=\rho_{+}-\rho_{-}$,etc.

Eliminating $\sigma$ from $(\ref{3.11})$, we obtain
\begin{equation}\label{3.12}
u-u_{-}=\pm\sqrt{\frac{\rho-\rho_{-}}{\rho\rho_{-}}\Big(A(\rho^{n}-\rho_{-}^{n})-B(\frac{1}{\rho^{\alpha}}-\frac{1}{\rho_{-}^{\alpha}})\Big)}.
\end{equation}

Using the Lax entropy condition, the 1-shock satisfies
\begin{equation}\label{3.13}
\sigma^{AB}<\lambda_{1}^{AB}(\rho_{-},u_{-}),\ \
\lambda_{1}^{AB}(\rho,u)<\sigma^{AB}<\lambda_{2}^{AB}(\rho,u),
\end{equation}
while the 1-shock satisfies
\begin{equation}\label{3.14}
\lambda_{1}^{AB}(\rho_{-},u_{-})<\sigma^{AB}<\lambda_{2}^{AB}(\rho_{-},u_{-}),\
\ \lambda_{2}^{AB}(\rho,u)<\sigma^{AB}.
\end{equation}

From the first equation in $(\ref{3.11})$, we have
\begin{equation}\label{3.15}
\sigma^{AB}=\frac{\rho
u-\rho_{-}u_{-}}{\rho-\rho_{-}}=u+\frac{\rho_{-}(
u-u_{-})}{\rho-\rho_{-}}=u_{-} +\frac{\rho(
u-u_{-})}{\rho-\rho_{-}}.
\end{equation}

Thus, by a simple calculation, $(\ref{3.13})$ is equivalent to
\begin{equation}\label{3.16}
-\rho\sqrt{An\rho^{n-1}+\frac{\alpha
B}{\rho^{\alpha+1}}}<\frac{\rho\rho_{-}(
u-u_{-})}{\rho-\rho_{-}}<-\rho_{-}\sqrt{An\rho_{-}^{n-1}+\frac{\alpha
B}{\rho_{-}^{\alpha+1}}},
\end{equation}
and $(\ref{3.14})$ is equivalent to
\begin{equation}\label{3.17}
\rho\sqrt{An\rho^{n-1}+\frac{\alpha
B}{\rho^{\alpha+1}}}<\frac{\rho\rho_{-}(
u-u_{-})}{\rho-\rho_{-}}<\rho_{-}\sqrt{An\rho_{-}^{n-1}+\frac{\alpha
B}{\rho_{-}^{\alpha+1}}}.
\end{equation}
$(\ref{3.16})$ and $(\ref{3.17})$ imply that $\rho>\rho_{-}$,
$u<u_{-}$ and $\rho<\rho_{-}$, $u<u_{-}$, respectively.

 Through the above analysis, for a given left
state $(\rho_-,u_-)$, the shock curves in the phase plane, which are
the sets of states that can be connected on the right by a 1-shock
or 2-shock, are as follows:

\
\begin{equation}\label{3.18}
S_1(\rho_-,u_-): \left\{\begin{array}{ll}
\sigma_1=\DF{\rho u-\rho_-u_-}{\rho-\rho_-},\\
u-u_-=-\sqrt{\frac{\rho-\rho_{-}}{\rho\rho_{-}}\Big(A(\rho^{n}-\rho_{-}^{n})-B(\frac{1}{\rho^{\alpha}}-\frac{1}{\rho_{-}^{\alpha}})\Big)},\
\ \rho>\rho_-,
\end{array}\right.
\end{equation}
and
\begin{equation}\label{3.19}
S_2(\rho_-,u_-): \left\{\begin{array}{ll}
\sigma_2=\DF{\rho u-\rho_-u_-}{\rho-\rho_-},\\
u-u_-=-\sqrt{\frac{\rho-\rho_{-}}{\rho\rho_{-}}\Big(A(\rho^{n}-\rho_{-}^{n})-B(\frac{1}{\rho^{\alpha}}-\frac{1}{\rho_{-}^{\alpha}})\Big)},\
\ \rho<\rho_-.
\end{array}\right.
\end{equation}

For the 1-shock wave, through differentiating $u$ respect to $\rho$
in the second equation in $(\ref{3.18})$, we get
\begin{equation}\label{3.20}
2(u-u_{-})u_{\rho}=\frac{1}{\rho^{2}}\Big(A(\rho^{n}-\rho_{-}^{n})-B(\frac{1}{\rho^{\alpha}}-\frac{1}{\rho_{-}^{\alpha}})\Big)+
\frac{\rho-\rho_{-}}{\rho\rho_{-}}(An\rho^{n-1}+\frac{\alpha
B}{\rho^{\alpha+1}})>0,
\end{equation}
which means that $u_{\rho}<0$ for the 1-shock wave and that the
1-shock wave curve is starlike with respect to $(\rho_-,u_-)$ in the
region $\rho>\rho_-$. Similarly, we can get  $u_{\rho}>0$ for the
2-shock wave and that the 2-shock wave curve is starlike with
respect to $(\rho_-,u_-)$ in the region $\rho<\rho_-$. In addition,
it is easy to check that
$\lim\limits_{\rho\rightarrow+\infty}u=-\infty$ for the 1-shock wave
and
 $\lim\limits_{\rho\rightarrow0}u=-\infty$ for the 2-shock wave.

Through the analysis above, for a given left state $(\rho_-,u_-)$,
the sets of states connected with $(\rho_-,u_-)$ on the right in the
phase plane consist of the 1-rarefaction wave curve
$R_1(\rho_-,u_-)$, the 2-rarefaction wave curve $R_2(\rho_-,u_-)$,
the 1-shock curve $S_1(\rho_-,u_-)$ and the 2-shock curve
$S_2(\rho_-,u_-)$. These curves divide the upper half plane into
four parts $R_{1}R_{2}(\rho_-,u_-)$, $R_{1}S_{2}(\rho_-,u_-)$,
$S_{1}R_{2}(\rho_-,u_-)$ and $S_{1}S_{2}(\rho_-,u_-)$. Now, we put
all of these curves together in the upper half plane ($\rho>0$,
$u\in R$) to obtain a picture as in Fig.2.

By the phase plane analysis method, it is easy to construct
Riemann solutions for any given right state $(\rho_+,u_+)$ as follows:\\
\ \ (1)  $(\rho_+,u_+)\in\rm{R_{1}R_{2}}(\rho_-,u_-):R_{1}+R_{2};$ \ \  (2)  $(\rho_+,u_+)\in\rm{R_{1}S_{2}}(\rho_-,u_-):R_{1}+S_{2};$\\
\ \ (3) $(\rho_+,u_+)\in\rm{S_{1}R_{2}}(\rho_-,u_-):S_{1}+R_{2};$ \
\ \
 (4) $(\rho_+,u_+)\in\rm{S_{1}S_{2}}(\rho_-,u_-):S_{1}+S_{2}.$


\unitlength 1.00mm 
\linethickness{0.4pt}
\ifx\plotpoint\undefined\newsavebox{\plotpoint}\fi 
\begin{picture}(98.71,55.84)(15,0)
\put(131.71,11.89){\vector(1,0){.07}}
\put(52.22,11.89){\line(1,0){79.49}}
\qbezier(57.57,16.35)(92.95,19.92)(111.08,49.84)
\qbezier(120.58,14.94)(78.31,20.94)(71.54,50.84)
\put(49.1,47.87){\vector(0,1){.07}}
\put(49.1,24.68){\line(0,1){23.19}}
\put(133.34,11.89){\makebox(0,0)[cc]{$u$}}
\put(49.1,49.59){\makebox(0,0)[cc]{$\rho$}}%
 \put(88,26.30){\circle*{.9}}%
 \put(96.71,26.02){\makebox(0,0)[cc]{$\scriptstyle(\rho_-,u_-)$}}
 \put(110.49,30.47){\makebox(0,0)[cc]{$R_{1}R_{2}(\rho_-,u_-)$}}
\put(87.3,42.96){\makebox(0,0)[cc]{$S_{1}R_{2}(\rho_-,u_-)$}}
\put(67.38,30.77){\makebox(0,0)[cc]{$S_{1}S_{2}(\rho_-,u_-)$}}
\put(84.92,16.2){\makebox(0,0)[cc]{$R_{1}S_{2}(\rho_-,u_-)$}}
\put(69.76,47.12){\makebox(0,0)[cc]{$S_1$}}
\put(111.99,46.08){\makebox(0,0)[cc]{$R_2$}}
\put(61.14,19.58){\makebox(0,0)[cc]{$S_2$}}
\put(110,19.43){\makebox(0,0)[cc]{$R_1$}}
 \put(94.41,3.99){\makebox(0,0)[cc]{ Fig.2}}
\end{picture}

\section{Formation of $\delta$-shocks and vacuum states as $A,B\rightarrow0$}

In this section, we will study the vanishing pressure limit process,
i.e.,$A,B\rightarrow0$. Since the two regions $S_{1}R_{2}
(\rho_-,u_-)$ and $R_{1}S_{2} (\rho_-,u_-)$ in the $(\rho,u)$ plane
have empty interior when $A,B\rightarrow0$, it suffices to analyze
the limit process for the two cases $(\rho_+,u_+)\in S_{1}S_{2}
(\rho_-,u_-)$ and $(\rho_+,u_+)\in R_{1}R_{2}(\rho_-,u_-)$.

Firstly, we analyze the formation of $\delta$-shocks in Riemann
solutions to the extended Chaplygin gas equations
$(\ref{1.1})$-$(\ref{1.2})$ with $(\ref{1.6})$ in the case
$(\rho_+,u_+)\in S_{1}S_{2}(\rho_-,u_-)$ as the pressure vanishes.

\subsection{Limit behavior of the Riemann solutions as $A,B\rightarrow0$}
When $(\rho_+,u_+)\in  S_{1}S_{2} (\rho_-,u_-)$, for fixed $A,B>0$,
let $(\rho_*^{AB},u_*^{AB})$ be the intermediate state in the sense
that $(\rho_-,u_-)$ and $(\rho_*^{AB},u_*^{AB})$ are connected by
1-shock $S_1$ with speed $\sigma_1^{AB}$, $(\rho_*^{AB},u_*^{AB})$
and $(\rho_+,u_+)$ are connected by 2-shock $S_2$ with speed
$\sigma_2^{AB}$. Then it follows
\begin{equation}\label{4.1}S_1:\ \
\left\{\begin{array}{ll}
\sigma_1^{AB}=\DF{\rho_{*}^{AB} u_{*}^{AB}-\rho_-u_-}{\rho_{*}^{AB}-\rho_-},\\
u_{*}^{AB}-u_-=-\sqrt{\frac{\rho_{*}^{AB}-\rho_{-}}{\rho_{*}^{AB}\rho_{-}}\Big(A((\rho_{*}^{AB})^{n}-\rho_{-}^{n})-
B(\frac{1}{(\rho_{*}^{AB})^{\alpha}}-\frac{1}{\rho_{-}^{\alpha}})\Big)},\
\ \rho_{*}^{AB}>\rho_-,
\end{array}\right.
\end{equation}

\begin{equation}\label{4.2}S_2:\ \
\left\{\begin{array}{ll}
\sigma_2^{AB}=\DF{\rho_{+} u_{+}-\rho_{*}^{AB}u_{*}^{AB}}{\rho_{+}-\rho_{*}^{AB}},\\
u_{+}-u_{*}^{AB}=-\sqrt{\frac{\rho_{+}-\rho_{*}^{AB}}{\rho_{+}\rho_{*}^{AB}}\Big(A(\rho_{+}^{n}-(\rho_{*}^{AB})^{n})-
B(\frac{1}{\rho_{+}^{\alpha}}-\frac{1}{(\rho_{*}^{AB})^{\alpha}})\Big)},\
\ \rho_{+}<\rho_{*}^{AB}.
\end{array}\right.
\end{equation}

In the following, we give some lemmas to show the limit behavior of
the Riemann solutions of system $(\ref{1.1})$-$(\ref{1.2})$ with
$(\ref{2.1})$ as $A,B\rightarrow0$.

\begin{Lemma}\label{lem:4.1}
$\lim\limits_{A,B\rightarrow0}\rho_*^{AB}=+\infty.$
\end{Lemma}

\noindent\textbf{Proof.} Eliminating $u_{*}^{AB}$ in the second
equation of $(\ref{4.1})$ and $(\ref{4.2})$ gives
\begin{eqnarray}
u_{+}-u_{-}=&-&\sqrt{\frac{\rho_{*}^{AB}-\rho_{-}}{\rho_{*}^{AB}\rho_{-}}
\Big(A((\rho_{*}^{AB})^{n}-\rho_{-}^{n})-B(\frac{1}{(\rho_{*}^{AB})^{\alpha}}-\frac{1}{\rho_{-}^{\alpha}})\Big)}\nonumber\\
&-&\sqrt{\frac{\rho_{+}-\rho_{*}^{AB}}{\rho_{+}\rho_{*}^{AB}}
\Big(A(\rho_{+}^{n}-(\rho_{*}^{AB})^{n})-B(\frac{1}{\rho_{+}^{\alpha}}-\frac{1}{(\rho_{*}^{AB})^{\alpha}})\Big)}.\label{4.3}
\end{eqnarray}
If
$\lim\limits_{A,B\rightarrow0}\rho_*^{AB}=M\in(\max\{\rho_{-},\rho_{+}\},+\infty)$,
then by taking the limit of $(\ref{4.3})$ as $A,B\rightarrow0$, we
obtain that $u_{+}-u_{-}=0$, which contradicts with $u_{+}<u_{-}$.
Therefore we must have
$\lim\limits_{A,B\rightarrow0}\rho_*^{AB}=+\infty.$

By Lemma \ref{lem:4.1}, from $(\ref{4.3})$ we immediately have the
following lemma.

\begin{Lemma}\label{lem:4.2}
$\lim\limits_{A,B\rightarrow0}A(\rho_*^{AB})^{n}=\DF{\rho_-\rho_+}{(\sqrt{\rho_-}+\sqrt{\rho_+})^2}(u_--u_+)^2$.
\end{Lemma}

\begin{Lemma}\label{lem:4.3}
\begin{equation}\label{4.4}
\lim\limits_{A,B\rightarrow0}u_*^{AB}=\lim\limits_{A,B\rightarrow0}\sigma_1^{AB}=\lim\limits_{A,B\rightarrow0}\sigma_2^{AB}=\sigma.
\end{equation}
\end{Lemma}

\noindent\textbf{Proof.} From the first equation of  $(\ref{4.1})$
and $(\ref{4.2})$ for $S_{1}$ and $S_{2}$, by Lemma \ref{lem:4.1},
we have
$$\lim\limits_{A,B\rightarrow0}\sigma_1^{AB}=
\lim\limits_{A,B\rightarrow0}\DF{\rho_{*}^{AB}
u_{*}^{AB}-\rho_-u_-}{\rho_{*}^{AB}-\rho_-}=
\lim\limits_{A,B\rightarrow0}\DF{
u_{*}^{AB}-\frac{\rho_-u_-}{\rho_{*}^{AB}}}{1-\frac{\rho_-}{\rho_{*}^{AB}}}=\lim\limits_{A,B\rightarrow0}u_*^{AB},$$

$$\lim\limits_{A,B\rightarrow0}\sigma_2^{AB}=
\lim\limits_{A,B\rightarrow0}\DF{\rho_+u_+-\rho_{*}^{AB}
u_{*}^{AB}}{\rho_+-\rho_{*}^{AB}}= \lim\limits_{A,B\rightarrow0}
\DF{
\frac{\rho_+u_+}{\rho_{*}^{AB}}-u_{*}^{AB}}{\frac{\rho_+}{\rho_{*}^{AB}}-1}=\lim\limits_{A,B\rightarrow0}u_*^{AB},$$
which immediately lead to
$\lim\limits_{A,B\rightarrow0}u_*^{AB}=\lim\limits_{A,B\rightarrow0}\sigma_1^{AB}=\lim\limits_{A,B\rightarrow0}\sigma_2^{AB}$.

From the second equation of $(\ref{4.1})$, by Lemma
\ref{lem:4.1}-\ref{lem:4.2}, we get

\begin{eqnarray}
\lim\limits_{A,B\rightarrow0}u_*^{AB}&=&u_{-}-\lim\limits_{A,B\rightarrow0}\sqrt{\frac{\rho_{*}^{AB}-\rho_{-}}{\rho_{*}^{AB}\rho_{-}}
\Big(A((\rho_{*}^{AB})^{n}-\rho_{-}^{n})-B(\frac{1}{(\rho_{*}^{AB})^{\alpha}}-\frac{1}{\rho_{-}^{\alpha}})\Big)}\nonumber\\
&=&u_{-}-\sqrt{\frac{1}{\rho_{-}}\DF{\rho_-\rho_+}{(\sqrt{\rho_-}+\sqrt{\rho_+})^2}(u_--u_+)^2}\nonumber\\
&=&u_{-}-\DF{\sqrt{\rho_+}}{\sqrt{\rho_-}+\sqrt{\rho_+}}(u_--u_+)\nonumber\\
&=&\frac{\sqrt{\rho_-}u_-+\sqrt{\rho_+}u_+}{\sqrt{\rho_-}+\sqrt{\rho_+}}=\sigma\nonumber.
\end{eqnarray}
The proof is completed.

\begin{Lemma}\label{lem:4.4}
\begin{equation}\label{4.5}
\lim\limits_{A,B\rightarrow0}\int_{\sigma_1^{AB}}^{\sigma_2^{AB}}\rho_*^{AB}d\xi=\sigma[\rho]-[\rho
u],
\end{equation}
\begin{equation}\label{4.6}
\lim\limits_{A,B\rightarrow0}\int_{\sigma_1^{AB}}^{\sigma_2^{AB}}\rho_*^{AB}u_*^{AB}d\xi=\sigma[\rho
u]-[\rho u^{2}].
\end{equation}

\end{Lemma}

\noindent\textbf{Proof.} The first equations of the Rankine-Hugoniot
condition $(\ref{3.11})$ for $S_{1}$ and $S_{2}$ read
\begin{equation}\label{4.7}
\left\{\begin{array}{ll}
\sigma_1^{AB}(\rho_*^{AB}-\rho_-)=\rho_{*}^{AB}
u_{*}^{AB}-\rho_-u_-,\\
\sigma_2^{AB}(\rho_+-\rho_{*}^{AB})=\rho_+u_+-\rho_{*}^{AB}
u_{*}^{AB},
\end{array}\right.
\end{equation}
from which we have
\begin{equation}\label{4.8}
\lim\limits_{A,B\rightarrow0}\rho_*^{AB}(\sigma_2^{AB}-\sigma_1^{AB})=\lim\limits_{A,B\rightarrow0}
(-\sigma_1^{AB}\rho_{-}+\sigma_2^{AB}\rho_{+}-\rho_+u_++\rho_-u_-)=\sigma[\rho]-[\rho
u].
\end{equation}

Similarly, from the second equations of the Rankine-Hugoniot
condition $(\ref{3.11})$ for $S_{1}$ and $S_{2}$
\begin{equation}\label{4.9}
\left\{\begin{array}{ll}
\sigma_1^{AB}(\rho_*^{AB}u_{*}^{AB}-\rho_-u_-)=\rho_{*}^{AB}
(u_{*}^{AB})^{2}-\rho_-u_-^{2}+A((\rho_{*}^{AB})^{\gamma}-\rho_-^{\gamma})-B(\DF{1}{(\rho_{*}^{AB})^{\alpha}}-\DF{1}{\rho_-^{\alpha}}),\\
\sigma_2^{AB}(\rho_+u_+-\rho_{*}^{AB}u_{*}^{AB})=\rho_+u_+^{2}-\rho_{*}^{AB}
(u_{*}^{AB})^{2}+A(\rho_+^{\gamma}-(\rho_{*}^{AB})^{\gamma})-B(\DF{1}{\rho_+^{\alpha}}-\DF{1}{(\rho_{*}^{AB})^{\alpha}}),
\end{array}\right.
\end{equation}

we obtain
\begin{eqnarray}
&\lim\limits_{A,B\rightarrow0}&\rho_*^{AB}u_{*}^{AB}(\sigma_2^{AB}-\sigma_1^{AB})\nonumber\\
&=&\lim\limits_{A,B\rightarrow0}(-\sigma_1^{AB}\rho_{-}u_{-}+\sigma_2^{AB}\rho_{+}u_{+}-
\rho_+u_+^{2}+\rho_-u_-^{2}-A(\rho_{+}^{\gamma}-\rho_{-}^{\gamma})+B(\DF{1}{\rho_+^{\alpha}}-\DF{1}{\rho_{-}^{\alpha}}))\nonumber\\
&=&\sigma[\rho u]-[\rho u^{2}].\label{4.10}
\end{eqnarray}
Thus, from $(\ref{4.8})$ and $(\ref{4.10})$ we immediately get
$(\ref{4.5})$ and $(\ref{4.6})$. The proof is finished.

\begin{Remark}
The above lemmas shows that, as $A,B\rightarrow0$, $S_{1}$ and
$S_{2}$ coincide, the intermediate density $\rho_{*}^{AB}$ becomes
singular, the velocities $\sigma_{1}^{AB}$, $\sigma_{2}^{AB}$ and
$u_{*}^{AB}$ for Riemann solutions of $(\ref{1.1})$-$(\ref{1.2})$
approach to $\sigma$, which are consistent with the velocity and the
density of the $\delta$-shock solution to the transport equations
$(\ref{1.3})$ with the same Riemann data $(\rho_\pm,u_\pm)$ in
Section 2.
\end{Remark}

\subsection{$\delta-$shocks and concentration}

Now we show the following theorem which is similar to Theorem 3.1 in
\cite{Chen-Liu1} and characterizes the vanishing pressure limit in
the case $(\rho_+,u_+)\in\rm{I\!V}(\rho_-,u_-)$ .

\begin{Theorem}\label{thm:4.1}
Let $u_->u_+$ and $(\rho_+,u_+)\in\rm{I\!V}(\rho_-,u_-)$. For any
fixed $A,B>0$, assume that $(\rho^{AB},u^{AB})$ is the two-shock
Riemann solution of $(\ref{1.1})$-$(\ref{1.2})$ with Riemann data
$(\rho_\pm,u_\pm)$ constructed in section 3. Then as
$A,B\rightarrow0$, $\rho^{AB}$ and $\rho^{AB}u^{AB}$ converge in the
sense of distributions, and the limit functions of $\rho^{AB}$ and
$\rho^{AB}u^{AB}$ are the sums of a step function and a
$\delta$-measure with weights
$$\frac{t}{\sqrt{1+\sigma^2}}(\sigma[\rho]-[\rho u])\ \ and\ \ \frac{t}{\sqrt{1+\sigma^2}}(\sigma[\rho u]-[\rho u^2]),$$
respectively, which form a $\delta$-shock solution of $(\ref{1.3})$
with the same Riemann data $(\rho_\pm,u_\pm)$.
\end{Theorem}

\noindent\textbf{Proof.} 1. Set $\xi=\frac{x}{t}$, for any fixed
$A,B>0$, the two-shock Riemann solution can be written as
\begin{equation}\nonumber
(\rho^{AB},u^{AB})(\xi)=\left\{\begin{array}{ll}
(\rho_-,u_-),\ \ \ \ \ \ -\infty<\xi<\sigma_1^{AB},\\
(\rho_*^{AB},u_*^{AB}),\ \ \ \sigma_1^{AB}<\xi<\sigma_2^{AB},\\
(\rho_+,u_+),\ \ \ \ \ \ \sigma_2^{AB}<\xi<\infty,
\end{array}\right.
\end{equation}
which satisfies the following weak formulations:
\begin{equation}\label{4.11}
-\int_{-\infty}^\infty(u^{AB}(\xi)-\xi)\rho^{AB}(\xi)\psi'(\xi)d\xi+\int_{-\infty}^\infty\rho^{AB}(\xi)\psi(\xi)d\xi=0,
\end{equation}
\begin{eqnarray}
-\int_{-\infty}^\infty(u^{AB}(\xi)-\xi)\rho^{AB}(\xi)u^{AB}(\xi)\psi'(\xi)d\xi+
\int_{-\infty}^\infty\rho^{AB}(\xi)u^{AB}(\xi)\psi(\xi)d\xi\nonumber\\
=\int_{-\infty}^\infty\Big(A(\rho^{AB}(\xi))^n-\DF{B}{(\rho^{AB}(\xi))^\alpha}\Big)\psi'(\xi)d\xi\label{4.12},
\end{eqnarray}
for any test function $\psi\in C_0^\infty(-\infty,\infty)$.

2. By using the weak formulation $(\ref{4.11})$, we can obtain the
limit of $\rho^{AB}$, which is denoted by the following identities:
\begin{equation}\label{4.13}
\lim\limits_{{A,B}\rightarrow0}\int_{-\infty}^\infty\Big(\rho^{AB}(\xi)-\rho_0(\xi-\sigma)\Big)\psi(\xi)d\xi=(\sigma[\rho]-[\rho
u])\psi(\sigma),
\end{equation}
for any test function $\psi\in C_0^\infty(-\infty,\infty)$, where
$$\rho_{0}(\xi)=\rho_{-}+[\rho]\chi(\xi),$$ and $\chi(\xi)$ is the characteristic function. Since the proof of $(\ref{4.13})$ is the same as step 2  in the
proof of Theorem 3.1 in \cite{Chen-Liu1}, we omit it.

3. Now we turn to justify the limit of $\rho^{AB} u^{AB}$ by using
the weak formulation $(\ref{4.12})$. The first integral on the left
hand side of $(\ref{4.12})$ can be decomposed into
\begin{equation}\label{4.14}
-\Big\{\int_{-\infty}^{\sigma_1^{AB}}+\int_{\sigma_1^{AB}}^{\sigma_2^{AB}}+\int_{\sigma_2^{AB}}^\infty\Big\}
(u^{AB}(\xi)-\xi)\rho^{AB}(\xi)u^{AB}(\xi)\psi'(\xi)d\xi.
\end{equation}

The sum of the first and last term of $(\ref{4.14})$ is
\begin{eqnarray}
&&-\int_{-\infty}^{\sigma_1^{AB}}(u_--\xi)\rho_-u_-\psi'(\xi)d\xi-\int_{\sigma_2^{AB}}^\infty(u_+-\xi)\rho_+u_+\psi'(\xi)d\xi\nonumber\\
&=&-\rho_-u_-^{2}\psi(\sigma_1^{AB})+\rho_+u_+^{2}\psi(\sigma_2^{AB})+\rho_-u_-
\sigma_1^{AB}\psi(\sigma_1^{AB})-\rho_+u_+\sigma_2^{AB}\psi(\sigma_2^{AB})\nonumber \\
&-&\int_{-\infty}^{\sigma_1^{AB}}\rho_-u_-\psi(\xi)d\xi-\int_{\sigma_2^{AB}}^\infty\rho_+u_+\psi(\xi)d\xi\nonumber,
\end{eqnarray}
which converges as $A,B\rightarrow0$ to
$$([\rho u^{2}]-\sigma[\rho u])\psi(\sigma)-\int_{-\infty}^\infty(\rho_0u_0)(\xi-\sigma)\psi(\xi)d\xi.$$

The second term of $(\ref{4.14})$ is
\begin{eqnarray}
&&-\rho_*^{AB}u_*^{AB}\int_{\sigma_1^{AB}}^{\sigma_2^{AB}}
(u_*^{AB}-\xi)\psi'(\xi)d\xi\nonumber\\
&=&-\rho_*^{AB}u_*^{AB}\Big((u_*^{AB}-\sigma_2^{AB})\psi(\sigma_2^{AB})-(u_*^{AB}-
\sigma_1^{AB})\psi(\sigma_1^{AB})+\int_{\sigma_1^{AB}}^{\sigma_2^{AB}}\psi(\xi)d\xi\Big)\nonumber\\
&=&-\rho_*^{AB}u_*^{AB}(\sigma_2^{AB}-\sigma_1^{AB})
\Big(u_*^{AB}\frac{\psi(\sigma_2^{AB})-\psi(\sigma_1^{AB})}{\sigma_2^{AB}-\sigma_1^{AB}}-
\frac{\sigma_2^{AB}\psi(\sigma_2^{AB})-\sigma_1^{AB}\psi(\sigma_1^{AB})}{\sigma_2^{AB}-\sigma_1^{AB}}\nonumber\\
&&+\frac{1}{\sigma_2^{AB}-\sigma_1^{AB}}\int_{\sigma_1^{AB}}^{\sigma_2^{AB}}\psi(\xi)d\xi
\Big)\nonumber,
\end{eqnarray}
which converges as $A,B\rightarrow0$ to
$$-(\sigma[\rho u-\rho u^{2}])\Big(\sigma\psi'(\sigma)-\sigma\psi'(\sigma)-\psi(\sigma)+\psi(\sigma)\Big)=0,$$
by Lemma \ref{lem:4.3}-\ref{lem:4.4}.

Now we compute the integral on the right hand side of
$(\ref{4.12})$, by Lemma \ref{lem:4.1}-\ref{lem:4.3}, we obtain
\begin{eqnarray}
&&\int_{-\infty}^\infty\Big(A(\rho^{AB}(\xi))^n-\DF{B}{(\rho^{AB}(\xi))^\alpha}\Big)\psi'(\xi)d\xi\nonumber\\
&=&\Big\{\int_{-\infty}^{\sigma_1^{AB}}+\int_{\sigma_1^{AB}}^{\sigma_2^{AB}}+\int_{\sigma_2^{AB}}^\infty\Big\}
\Big(A(\rho^{AB}(\xi))^n-\DF{B}{(\rho^{AB}(\xi))^\alpha}\Big)\psi'(\xi)d\xi\nonumber\\
&=&\int_{-\infty}^{\sigma_1^{AB}}\Big(A\rho_{-}^n-\DF{B}{\rho_{-}^\alpha}\Big)\psi'(\xi)d\xi+
\int_{\sigma_1^{AB}}^{\sigma_2^{AB}}\Big(A(\rho_{*}^{AB})^n-\DF{B}{(\rho_{*}^{AB})^\alpha}\Big)\psi'(\xi)d\xi\nonumber\\
&&+\int_{\sigma_2^{AB}}^\infty\Big(A\rho_{+}^n-\DF{B}{\rho_{+}^\alpha}\Big)\psi'(\xi)d\xi\nonumber\\
&=&\Big(A\rho_{-}^n-\DF{B}{\rho_{-}^\alpha}\Big)\psi(\sigma_1^{AB})
-\Big(A\rho_{+}^n-\DF{B}{\rho_{+}^\alpha}\Big)\psi(\sigma_2^{AB})+
\Big(A(\rho_{*}^{AB})^n-\DF{B}{(\rho_{*}^{AB})^\alpha}\Big)(\psi(\sigma_2^{AB})-\psi(\sigma_1^{AB})),
\nonumber
\end{eqnarray}
which converge to 0  as $A,B\rightarrow0$.

Then, the integral identity $(\ref{4.12})$ yields
\begin{equation}\label{4.15}
\lim\limits_{A,B\rightarrow0}\int_{-\infty}^\infty\Big((\rho^{AB}
u^{AB})(\xi)-(\rho_0u_0)(\xi-\sigma)\Big)\psi(\xi)d\xi=(\sigma[\rho
u]-[\rho u^{2}])\psi(\sigma),
\end{equation}
for any test function $\psi\in C_0^\infty(-\infty,\infty)$.

4. Finally, we are in the position to study the limits of
$\rho^{AB}$ and $\rho^{AB}u^{AB}$ by tracking the time-dependence of
the weights of the $\delta$-measure as $A,B\rightarrow0$.

Let $\phi(x,t)\in  C_0^\infty((-\infty,\infty)\times[0,\infty))$ be
a smooth test function and $\tilde{\phi}(\xi,t)=\phi(\xi t,t)$. Then
we have
$$\lim\limits_{A,B\rightarrow0}\int_0^\infty\int_{-\infty}^\infty\rho^{AB}(\frac{x}{t})\phi(x,t)dxdt
=\lim\limits_{A,B\rightarrow0}\int_0^\infty
t\Big(\int_{-\infty}^\infty\rho^{AB}(\xi)\tilde{\phi}(\xi,t)d\xi\Big)dt.$$
On the other hand, from $(\ref{4.13})$, we have
\begin{eqnarray*}
\lim\limits_{A,B\rightarrow0}\int_{-\infty}^\infty\rho^{AB}(\xi)\tilde{\phi}(\xi,t)d\xi
&=&\int_{-\infty}^\infty\rho_0(\xi-\sigma)\tilde{\phi}(\xi,t)d\xi+(\sigma[\rho ]-[\rho u])\tilde{\phi}(\sigma,t)\\
&=&\frac{1}{t}\int_{-\infty}^\infty\rho_0(x-\sigma
t)\phi(x,t)dx+(\sigma[\rho ]-[\rho u])\phi(\sigma t,t).
\end{eqnarray*}
Combining the two relations above yields
\begin{eqnarray*}
\lim\limits_{A,B\rightarrow0}\int_0^\infty\int_{-\infty}^\infty\rho^{AB}(\frac{x}{t})\phi(x,t)dxdt=\int_0^\infty\int_{-\infty}^\infty\rho_0(x-\sigma
t)\phi(x,t)dxdt+\int_0^\infty t(\sigma[\rho ]-[\rho u])\phi(\sigma
t,t)dt.
\end{eqnarray*}
The last term, by definition, equals to
$$\langle w_1(\cdot)\delta_S,\phi(\cdot,\cdot)\rangle,$$
with
$$w_1(t)=\frac{t}{\sqrt{1+\sigma^2}}(\sigma[\rho]-[\rho u]).$$

Similarly, from $(\ref{4.15})$ we can show that
\begin{eqnarray*}
\lim\limits_{A,B\rightarrow0}\int_0^\infty\int_{-\infty}^\infty(\rho^{AB}
u^{AB})(\frac{x}{t})\phi(x,t)dxdt=\int_0^\infty\int_{-\infty}^\infty(\rho_0u_0)(x-\sigma
t)\phi(x,t)dxdt+\langle w_2(\cdot)\delta_S,\phi(\cdot,\cdot)\rangle,
\end{eqnarray*}
with
$$w_2(t)=\frac{t}{\sqrt{1+\sigma^2}}(\sigma[\rho u]-[\rho u^2]).$$
Then we complete the proof of Theorem \ref{thm:4.1}.

\subsection{Formation of vacuum states} In this subsection, we show the
formation of vacuum states in the Riemann solutions to
$(\ref{1.1})$-$(\ref{1.2})$ with $(\ref{1.6})$ in the case
$(\rho_+,u_+)\in R_{1}R_{2}(\rho_-,u_-)$ with $u_-<u_+$ and
$\rho_\pm>0$ as the pressure vanishes.

At this monent, for fixed $A,B>0$, let $(\rho_*^{AB},u_*^{AB})$ be
the intermediate state in the sense that $(\rho_-,u_-)$ and
$(\rho_*^{AB},u_*^{AB})$ are connected by 1-rarefaction wave $R_1$
with speed $\lambda_1^{AB}$, $(\rho_*^{AB},u_*^{AB})$and
$(\rho_+,u_+)$ are connected by 2-rarefaction wave $R_2$ with speed
$\lambda_2^{AB}$. Then it follows
\begin{equation}\label{4.16}R_1:\ \
\left\{\begin{array}{ll} \xi=\lambda_1^{AB}=u-\sqrt{An\rho^{n-1}+\frac{\alpha B}{\rho^{\alpha+1}}},\\
u-u_-=-\int_{\rho_{-}}^{\rho}\frac{\sqrt{An\rho^{n-1}+\frac{\alpha
B}{\rho^{\alpha+1}}}}{\rho}d\rho,\ \
\rho_{*}^{AB}\leq\rho\leq\rho_{-}.
\end{array}\right.
\end{equation}

\begin{equation}\label{4.17}R_2:\ \
\left\{\begin{array}{ll} \xi=\lambda_2^{AB}=u+\sqrt{An\rho^{n-1}+\frac{\alpha B}{\rho^{\alpha+1}}},\\
u_+-u=\int_{\rho}^{\rho_{+}}\frac{\sqrt{An\rho^{n-1}+\frac{\alpha
B}{\rho^{\alpha+1}}}}{\rho}d\rho,\ \
\rho_{*}^{AB}\leq\rho\leq\rho_{+}.
\end{array}\right.
\end{equation}

Now, from the second equations of $(\ref{4.16})$ and $(\ref{4.17})$,
using the following integral identity
\begin{eqnarray}
&&\int^{\rho_{-}}_{\rho}\frac{\sqrt{An\rho_{-}^{n-1}+\frac{\alpha
B}{\rho^{\alpha+1}}}}{\rho}d\rho\nonumber\\
&=&\frac{2}{\alpha+1}\Big(-\sqrt{An\rho_{-}^{n-1}+\frac{\alpha
B}{\rho^{\alpha+1}}} +\sqrt{An\rho_{-}^{n-1}}
\ln(\sqrt{An\rho_{-}^{n-1}\rho^{\alpha+1}+\alpha
B}+\sqrt{An\rho_{-}^{n-1}\rho^{\alpha+1}})\Big)\Big
|_{\rho}^{\rho_{-}}\nonumber,
\end{eqnarray}

 it follows that the
intermediate state $(\rho_*^{AB},u_*^{AB})$ satisfies

\begin{eqnarray}
&&u_+-u_-\nonumber\\
&=&\int^{\rho_{-}}_{\rho_*^{AB}}\frac{\sqrt{An\rho^{n-1}+\frac{\alpha
B}{\rho^{\alpha+1}}}}{\rho}d\rho+\int_{\rho_*^{AB}}^{\rho_{+}}\frac{\sqrt{An\rho^{n-1}+\frac{\alpha
B}{\rho^{\alpha+1}}}}{\rho}d\rho\nonumber\\
&\leq&\int^{\rho_{-}}_{\rho_*^{AB}}\frac{\sqrt{An\rho_{-}^{n-1}+\frac{\alpha
B}{\rho^{\alpha+1}}}}{\rho}d\rho+\int_{\rho_*^{AB}}^{\rho_{+}}\frac{\sqrt{An\rho_{+}^{n-1}+\frac{\alpha
B}{\rho^{\alpha+1}}}}{\rho}d\rho\nonumber\\
&=&\frac{2}{\alpha+1}\Big(-\sqrt{An\rho_{-}^{n-1}+\frac{\alpha
B}{\rho_{-}^{\alpha+1}}} +\sqrt{An\rho_{-}^{n-1}}
\ln(\sqrt{An\rho_{-}^{n-1}\rho_{-}^{\alpha+1}+\alpha
B}+\sqrt{An\rho_{-}^{n-1}\rho_{-}^{\alpha+1}})\nonumber\\
&&+\sqrt{An\rho_{-}^{n-1}+\frac{\alpha B}{(\rho_*^{AB})^{\alpha+1}}}
-\sqrt{An\rho_{-}^{n-1}}
\ln(\sqrt{An\rho_{-}^{n-1}(\rho_*^{AB})^{\alpha+1}+\alpha
B}+\sqrt{An\rho_{-}^{n-1}(\rho_*^{AB})^{\alpha+1}})\nonumber\\
&&-\sqrt{An\rho_{+}^{n-1}+\frac{\alpha B}{\rho_{+}^{\alpha+1}}}
+\sqrt{An\rho_{+}^{n-1}}
\ln(\sqrt{An\rho_{+}^{n-1}\rho_{+}^{\alpha+1}+\alpha
B}+\sqrt{An\rho_{+}^{n-1}\rho_{+}^{\alpha+1}})\nonumber\\
&&+\sqrt{An\rho_{+}^{n-1}+\frac{\alpha B}{(\rho_*^{AB})^{\alpha+1}}}
-\sqrt{An\rho_{+}^{n-1}}
\ln(\sqrt{An\rho_{+}^{n-1}(\rho_*^{AB})^{\alpha+1}+\alpha
B}+\sqrt{An\rho_{+}^{n-1}(\rho_*^{AB})^{\alpha+1}})\Big
),\nonumber\\\label{4.18}
\end{eqnarray}
which implies the following result.

\begin{Theorem}\label{thm:4.2}
Let $u_-<u_+$ and $(\rho_+,u_+)\in\rm{I}(\rho_-,u_-)$. For any fixed
$A,B>0$, assume that $(\rho^{AB},u^{AB})$ is the two-rarefaction
wave Riemann solution of $(\ref{1.1})$-$(\ref{1.2})$ with Riemann
data $(\rho_\pm,u_\pm)$ constructed in section 3. Then as
$A,B\rightarrow0$, the limit of the Riemann solution
 $(\rho^{AB}, u^{AB})$ is two contact
discontinuities connecting the constant states $(\rho_\pm,u_\pm)$
and the intermediate vacuum state as follows:
\begin{equation}\nonumber
(\rho,u)(\xi)=\left\{\begin{array}{ll}
(\rho_-,u_-),\ \ \ \ \ \ -\infty<\xi\leq u_-,\\
(0,\xi),\ \ \ \ \ \ \ \ \ \ u_-\leq\xi\leq u_+,\\
(\rho_+,u_+),\ \ \ \ \ \ u_+\leq\xi<\infty,
\end{array}\right.
\end{equation}
which is exactly the Riemann solution to the transport equations
$(\ref{1.3})$ with the same Riemann data $(\rho_\pm,u_\pm)$.
\end{Theorem}

Indeed, if
$\lim\limits_{A,B\rightarrow0}\rho_*^{AB}=K\in(0,\min\{\rho_{-},\rho_{+}\})$,
then $(\ref{4.18})$ leads to $u_+-u_-=0$, which contradicts with
$u_-<u_+$. Thus $\lim\limits_{A,B\rightarrow0}\rho_*^{AB}=0$, which
just means vacuum occurs. Moreover, as $A,B\rightarrow0$, one can
directly derive from $(\ref{4.16})$ and $(\ref{4.17})$ that
$\lambda_{1}^{AB},\ \lambda_{2}^{AB}\rightarrow u$ and two
rarefaction waves $R_{1}$ and $R_{2}$ tend to two contact
discontinuities $\xi=\frac{x}{t}=u_\pm$, respectively. These reach
the desired conclusion.

\section{Formation of $\delta$-shocks and two-rarefaction wave as $A\rightarrow0$}
In this section, we discuss the limit behaviors of Riemann solutions
of $(\ref{1.1})$-$(\ref{1.2})$ with $(\ref{1.6})$ as the pressure
approaches the generalized Chaplygin gas pressure, i.e.,
$A\rightarrow0$.

From Section 2 and 3, we can easily check that, as $A\rightarrow0$,
the backward (forward) rarefaction wave curve $R_{1} (R_{2})$ of
$(\ref{1.1})$-$(\ref{1.2})$ tends to the backward (forward)
rarefaction wave curve $\overleftarrow{R }(\overrightarrow{R})$ of
$(\ref{1.5})$, and the backward (forward) shock wave curve $S_{1}
(S_{2})$ of $(\ref{1.1})$-$(\ref{1.2})$ tends to the backward
(forward) rarefaction wave curve $\overleftarrow{S
}(\overrightarrow{S})$ of $(\ref{1.5})$ when $0<\alpha<1$, while the
backward (forward) rarefaction wave curve $R_{1} (R_{2})$ of
$(\ref{1.1})$-$(\ref{1.2})$ tends to the backward (forward) contact
discontinuity curve of $(\ref{1.5})$, and the backward (forward)
shock wave curve $S_{1} (S_{2})$ of $(\ref{1.1})$-$(\ref{1.2})$
tends to the backward (forward) contact discontinuity curve of
$(\ref{1.5})$ when $\alpha=1$ (see Fig.3).

\unitlength 1mm 
\linethickness{0.4pt}
\ifx\plotpoint\undefined\newsavebox{\plotpoint}\fi 
\begin{picture}(115,72)(-8,0)
\put(110.75,12){\vector(1,0){.07}} \put(13,12){\line(1,0){97.75}}
\put(9.25,72){\vector(0,1){.07}} \put(9.25,15){\line(0,1){57}}

\bezier{50}(52,68)(52,50)(52,12) \bezier{50}(90,68)(90,50)(90,12)
 \qbezier(100,60) (78.75,16.5)(28,18.5)
 \qbezier(40,55) (69,14.5)(107,19)

 \qbezier(89.5,67.5)(78.13,15)(39.25,13.5)
\qbezier(52.5,67.5)(60.13,13.5)(103.25,14.5)
\qbezier(50.75,67)(43,14.5)(14.25,13) \tiny{
\put(113,11.5){\makebox(0,0)[cc]{$u$}}
\put(7,70.25){\makebox(0,0)[cc]{$\rho$}}
\put(51.5,8){\makebox(0,0)[cc]{$u_--\sqrt{
B}\rho_-^{-\frac{\alpha+1}{2}}$}}
\put(89.25,8){\makebox(0,0)[cc]{$u_-+\sqrt{
B}\rho_-^{-\frac{\alpha+1}{2}}$}}
\put(94.75,33.25){\makebox(0,0)[cc]{I}}
\put(71.75,16.75){\makebox(0,0)[cc]{I\!I}}
\put(70.25,44.75){\makebox(0,0)[cc]{I\!I\!I}}
\put(70,26.5){\circle*{.6}} \put(37.5,27){\circle*{.6}}
\put(23,27){\makebox(0,0)[cc]{$(\rho_-,u_--2\sqrt{
B}\rho_-^{-\frac{\alpha+1}{2}})$}}
\put(81,26.5){\makebox(0,0)[cc]{$(\rho_-,u_-)$}}
\put(48.25,30.25){\makebox(0,0)[cc]{I\!V}}
\put(25.75,37.5){\makebox(0,0)[cc]{V}}
\put(58,58.5){\makebox(0,0)[cc]{$\overleftarrow{S}$}}
\put(39,53.5){\makebox(0,0)[cc]{$S_{1}$}}
\put(39,20){\makebox(0,0)[cc]{$S_{2}$}}
\put(100,53.5){\makebox(0,0)[cc]{$R_{2}$}}
\put(100,20){\makebox(0,0)[cc]{$R_{1}$}}
\put(84,58.25){\makebox(0,0)[cc]{$\overrightarrow{R}$}}
\put(95,13.5){\makebox(0,0)[cc]{$\overleftarrow{R}$}}
\put(50,13.5){\makebox(0,0)[cc]{$\overrightarrow{S}$}}
\put(40,40){\makebox(0,0)[cc]{$S_{\delta}$}}
\put(52,2){\makebox(0,0)[cc]{Fig.1}}}
\end{picture}


\subsection{Formation of $\delta$-shocks}
In this subsection, we study the formation of the delta shock waves
in the limit as $A\rightarrow0$ of solutions of
$(\ref{1.1})$-$(\ref{1.2})$ with $(\ref{1.6})$ in the case
$(\rho_+,u_+)\in \rm{V} (\rho_-,u_-)$, i.e.,  $u_++\sqrt{
B}\rho_+^{-\frac{\alpha+1}{2}} \leq u_--\sqrt{
B}\rho_-^{-\frac{\alpha+1}{2}}$.

\begin{Lemma}\label{lem:5.1}
When $(\rho_+,u_+)\in \rm{V} (\rho_-,u_-)$, there exists a positive
parameter $A_0$ such that $(\rho_+,u_+)\in S_{1}S_{2} (\rho_-,u_-)$
when $0<A<A_{0}$.
\end{Lemma}

\noindent\textbf{Proof.} From $(\rho_+,u_+)\in \rm{V} (\rho_-,u_-)$,
we have

\begin{equation}\label{5.1}
u_++\sqrt{ B}\rho_+^{-\frac{\alpha+1}{2}} \leq u_--\sqrt{
B}\rho_-^{-\frac{\alpha+1}{2}},
\end{equation}
then
\begin{eqnarray}
( u_--u_+)^{2}&\geq&\big(\sqrt{ B}\rho_+^{-\frac{\alpha+1}{2}}
+\sqrt{
B}\rho_-^{-\frac{\alpha+1}{2}}\big)^{2}\nonumber\\
&=&B(\rho_+^{-\alpha-1}+\rho_-^{-\alpha-1}+2\rho_{+}^{-\frac{\alpha+1}{2}}\rho_{-}^{-\frac{\alpha+1}{2}})\nonumber\\
&>&B(\rho_+^{-\alpha-1}+\rho_-^{-\alpha-1}-\rho_{+}^{-1}\rho_{-}^{-\alpha}-\rho_{-}^{-1}\rho_{+}^{-\alpha})\nonumber\\
&=&B(\frac{1}{\rho_{+}}-\frac{1}{\rho_{-}})(\frac{1}{\rho_{+}^{\alpha}}-\frac{1}{\rho_{-}^{\alpha}})\label{5.2}.
\end{eqnarray}

All the states $(\rho,u)$ connected with $(\rho_-,u_-)$ by a
backward shock wave $S_{1}$ or a forward shock wave $S_{2}$ satisfy
\begin{equation}\label{5.3}
u-u_-=-\sqrt{\frac{\rho-\rho_{-}}{\rho\rho_{-}}\Big(A(\rho^{n}-\rho_{-}^{n})-B(\frac{1}{\rho^{\alpha}}-\frac{1}{\rho_{-}^{\alpha}})\Big)},\
\ \rho>\rho_-,
\end{equation}
or
\begin{equation}\label{5.4}
u-u_-=-\sqrt{\frac{\rho-\rho_{-}}{\rho\rho_{-}}\Big(A(\rho^{n}-\rho_{-}^{n})-B(\frac{1}{\rho^{\alpha}}-\frac{1}{\rho_{-}^{\alpha}})\Big)},\
\ \rho<\rho_-.
\end{equation}

When $\rho_{+}=\rho_{-}$, the conclusion is obviously true.  When
$\rho_{+}\neq\rho_{-}$, by taking
\begin{equation}\label{5.5}
(u_+-u_-)^{2}=\frac{\rho_+-\rho_{-}}{\rho_+\rho_{-}}\Big(A_{0}(\rho_+^{n}-\rho_{-}^{n})-B(\frac{1}{\rho_{+}^{\alpha}}-\frac{1}{\rho_{-}^{\alpha}})\Big),
\end{equation}
we have
\begin{equation}\label{5.6}
A_{0}=\frac{\rho_{+}\rho_{-}}{(\rho_{+}-\rho_{-})(\rho_{+}^{n}-\rho_{-}^{n})}
\Big((u_+-u_-)^{2}-B(\frac{1}{\rho_{+}}-\frac{1}{\rho_{-}})(\frac{1}{\rho_{+}^{\alpha}}-\frac{1}{\rho_{-}^{\alpha}})\Big),
\end{equation}
which together with $(\ref{5.2})$ gives the conclusion. The proof is
completed.

When $0<A<A_{0}$, the Riemann solution of
$(\ref{1.1})$-$(\ref{1.2})$ with $(\ref{1.6})$ includes a backward
shock wave $S_{1}$ and a forward shock wave $S_{2}$ with the
intermediate state $(\rho_{*}^{A},u_{*}^{A})$ besides two constant
states $(\rho_{\pm},u_{\pm})$. We then have

\begin{equation}\label{5.7}S_1:\ \
\left\{\begin{array}{ll}
\sigma_1^{AB}=\DF{\rho_{*}^{A} u_{*}^{A}-\rho_-u_-}{\rho_{*}^{A}-\rho_-},\\
u_{*}^{A}-u_-=-\sqrt{\frac{\rho_{*}^{A}-\rho_{-}}{\rho_{*}^{A}\rho_{-}}\Big(A((\rho_{*}^{A})^{n}-\rho_{-}^{n})-
B(\frac{1}{(\rho_{*}^{A})^{\alpha}}-\frac{1}{\rho_{-}^{\alpha}})\Big)},\
\ \rho_{*}^{A}>\rho_-,
\end{array}\right.
\end{equation}
and
\begin{equation}\label{5.8}S_2:\ \
\left\{\begin{array}{ll}
\sigma_2^{AB}=\DF{\rho_{+} u_{+}-\rho_{*}^{A}u_{*}^{A}}{\rho_{+}-\rho_{*}^{A}},\\
u_{+}-u_{*}^{A}=-\sqrt{\frac{\rho_{+}-\rho_{*}^{A}}{\rho_{+}\rho_{*}^{A}}\Big(A(\rho_{+}^{n}-(\rho_{*}^{A})^{n})-
B(\frac{1}{\rho_{+}^{\alpha}}-\frac{1}{(\rho_{*}^{A})^{\alpha}})\Big)},\
\ \rho_{+}<\rho_{*}^{A}.
\end{array}\right.
\end{equation}
Here $\sigma_1^{A}$ and $\sigma_2^{A}$ are the propagation speed of
$S_{1}$ and $S_{2}$, respectively. Similar to that in Section 4,
 in the following, we give some lemmas to show the limit behavior
of the Riemann solutions of system $(\ref{1.1})$-$(\ref{1.2})$ with
$(\ref{1.6})$ as $A\rightarrow0$.

\begin{Lemma}\label{lem:5.2}
$\lim\limits_{A\rightarrow0}\rho_*^{A}=+\infty.$
\end{Lemma}

\noindent\textbf{Proof.} Eliminating $u_{*}^{AB}$ in the second
equation of $(\ref{5.7})$ and $(\ref{5.8})$ gives
\begin{eqnarray}
u_{-}-u_{+}=&&\sqrt{\frac{\rho_{*}^{A}-\rho_{-}}{\rho_{*}^{A}\rho_{-}}
\Big(A((\rho_{*}^{A})^{n}-\rho_{-}^{n})-B(\frac{1}{(\rho_{*}^{A})^{\alpha}}-\frac{1}{\rho_{-}^{\alpha}})\Big)}\nonumber\\
&+&\sqrt{\frac{\rho_{+}-\rho_{*}^{A}}{\rho_{+}\rho_{*}^{A}}
\Big(A(\rho_{+}^{n}-(\rho_{*}^{A})^{n})-B(\frac{1}{\rho_{+}^{\alpha}}-\frac{1}{(\rho_{*}^{A})^{\alpha}})\Big)}.\label{5.9}
\end{eqnarray}
If
$\lim\limits_{A\rightarrow0}\rho_*^{A}=K\in(\max\{\rho_{-},\rho_{+}\},+\infty)$,
then by taking the limit of $(\ref{5.9})$ as $A\rightarrow0$, we
obtain that
\begin{eqnarray}
u_{-}-u_{+}&=&\sqrt{B}\Big(\sqrt{(\frac{1}{\rho_{-}}-\frac{1}{K})(\frac{1}{\rho_{-}^{\alpha}}-\frac{1}{K^{\alpha}})}+
\sqrt{(\frac{1}{\rho_{+}}-\frac{1}{K})(\frac{1}{\rho_{+}^{\alpha}}-\frac{1}{K^{\alpha}})}\Big)\nonumber\\
&<&\sqrt{B}\Big(\sqrt{\frac{1}{\rho_{-}}\frac{1}{\rho_{-}^{\alpha}}}+
\sqrt{\frac{1}{\rho_{+}}\frac{1}{\rho_{+}^{\alpha}}}\Big)\nonumber\\
&=&\sqrt{ B}\rho_+^{-\frac{\alpha+1}{2}} +\sqrt{
B}\rho_-^{-\frac{\alpha+1}{2}}\label{5.10}
\end{eqnarray}
which contradicts with $(\ref{5.1})$. Therefore we must have
$\lim\limits_{A\rightarrow0}\rho_*^{A}=+\infty.$

By Lemma \ref{lem:5.2}, from $(\ref{5.9})$ we immediately have the
following lemma.

\begin{Lemma}\label{lem:5.3}
$\lim\limits_{A\rightarrow0}A(\rho_*^{A})^{n}<\rho_-(u_--u_+)^2$.
\end{Lemma}

\begin{Lemma}\label{lem:5.4}
Let $\lim\limits_{A\rightarrow0}u_*^{A}=\widehat{\sigma^{B}}$, then
\begin{equation}\label{5.11}
\lim\limits_{A\rightarrow0}u_*^{A}=\lim\limits_{A\rightarrow0}\sigma_1^{A}
=\lim\limits_{A\rightarrow0}\sigma_2^{A}=\widehat{\sigma^{B}}\in\Big(u_++\sqrt{
\alpha B}\rho_+^{-\frac{\alpha+1}{2}}, u_--\sqrt{\alpha
B}\rho_-^{-\frac{\alpha+1}{2}}\Big).
\end{equation}
\end{Lemma}

\noindent\textbf{Proof.} From the second equation of $(\ref{5.7})$
for $S_{1}$, by Lemma \ref{lem:4.2} and \ref{lem:4.3}, we have
\begin{eqnarray}
\lim\limits_{A\rightarrow0}u_*^{A}&=&u_{-}-\lim\limits_{A\rightarrow0}\sqrt{\frac{\rho_{*}^{A}-\rho_{-}}{\rho_{*}^{A}\rho_{-}}
\Big(A((\rho_{*}^{A})^{n}-\rho_{-}^{n})-B(\frac{1}{(\rho_{*}^{A})^{\alpha}}-\frac{1}{\rho_{-}^{\alpha}})\Big)}\nonumber\\
&=&u_{-}-\sqrt{\frac{1}{\rho_{-}}\Big(\lim\limits_{A\rightarrow0}A(\rho_*^{A})^{n}+\frac{B}{\rho_{-}^{\alpha}}\Big)}\nonumber\\
&<& u_--\sqrt{\alpha
B}\rho_-^{-\frac{\alpha+1}{2}}\nonumber\\\label{5.12}.
\end{eqnarray}

Similarly, from the second equation of $(\ref{5.8})$ for $S_{2}$, we
have
\begin{eqnarray}
\lim\limits_{A\rightarrow0}u_*^{A}&=&u_{+}+\lim\limits_{A\rightarrow0}\sqrt{\frac{\rho_{+}-\rho_{*}^{A}}{\rho_{+}\rho_{*}^{A}}\Big(A(\rho_{+}^{n}-(\rho_{*}^{A})^{n})-
B(\frac{1}{\rho_{+}^{\alpha}}-\frac{1}{(\rho_{*}^{A})^{\alpha}})\Big)}\nonumber\\
&=&u_{+}+\sqrt{\frac{1}{\rho_{+}}\Big(\lim\limits_{A\rightarrow0}A(\rho_*^{A})^{n}+\frac{B}{\rho_{+}^{\alpha}}\Big)}\nonumber\\
&>&u_++\sqrt{ \alpha B}\rho_+^{-\frac{\alpha+1}{2}}
\nonumber\\\label{5.13}.
\end{eqnarray}
Furthermore, similar to the analysis as Lemma \ref{lem:4.3}, we can
obtain
$\lim\limits_{A\rightarrow0}u_*^{A}=\lim\limits_{A\rightarrow0}\sigma_1^{A}=\lim\limits_{A\rightarrow0}\sigma_2^{A}=\widehat{\sigma^{B}}$.
The proof is complete.

Similar to  Lemma \ref{lem:4.4}, we have the following lemma.
\begin{Lemma}\label{lem:5.5}
\begin{equation}\label{5.14}
\lim\limits_{A\rightarrow0}\int_{\sigma_1^{A}}^{\sigma_2^{A}}\rho_*^{A}d\xi=\sigma^{B}[\rho]-[\rho
u],
\end{equation}
\begin{equation}\label{5.15}
\lim\limits_{A\rightarrow0}\int_{\sigma_1^{A}}^{\sigma_2^{A}}\rho_*^{A}u_*^{A}d\xi=\sigma^{B}[\rho
u]-[\rho u^{2}-\frac{B}{\rho^{\alpha}}].
\end{equation}

\end{Lemma}

\begin{Lemma}\label{lem:5.6}
For $\widehat{\sigma^{B}}$ mentioned in Lemma \ref{lem:5.4},
\begin{equation}\label{5.16}
\widehat{\sigma^{B}}=\sigma^{B}=\DF{\rho_+ u_+-\rho_-u_-+
\big\{\rho_+\rho_-\big((u_+-u_-)^2-(\frac{1}{\rho_+}-\frac{1}{\rho_-})
(\frac{B}{\rho_+^\alpha}-\frac{B}{\rho_-^\alpha})\big)\big\}^{\frac{1}{2}}}{\rho_+-\rho_-},
\end{equation}
as $\rho_+\neq\rho_-$, and

\begin{equation}\label{5.17}
\widehat{\sigma^{B}}=\sigma^{B}=\frac{u_++u_-}{2}
\end{equation}
as $\rho_+=\rho_-$.

\end{Lemma}

\noindent\textbf{Proof.} Let
$\lim\limits_{A\rightarrow0}A(\rho_*^{A})^{n}=L$, by Lemma
\ref{lem:5.4}, from $(\ref{5.12})$ and $(\ref{5.13})$ we have
$$\lim\limits_{A\rightarrow0}u_*^{A}=
u_{-}-\sqrt{\frac{1}{\rho_{-}}\Big(L+\frac{B}{\rho_{-}^{\alpha}}\Big)}=
u_{+}+\sqrt{\frac{1}{\rho_{+}}\Big(L+\frac{B}{\rho_{+}^{\alpha}}\Big)}=\widehat{\sigma^{B}},$$
which leas to
\begin{equation}\label{5.18}
L+\frac{B}{\rho_{+}^{\alpha}}=\rho_{-}(u_--\widehat{\sigma^{B}})^{2},
\end{equation}
\begin{equation}\label{5.19}
L+\frac{B}{\rho_{-}^{\alpha}}=\rho_{+}(u_+-\widehat{\sigma^{B}})^{2}.
\end{equation}
Eliminating $L$ from $(\ref{5.18})$ and $(\ref{5.19})$, we have

\begin{equation}\label{5.20}
(\rho_+-\rho_-)(\widehat{\sigma^{B}})^{2}-2(\rho_+
u_+-\rho_-u_-)\widehat{\sigma^{B}}+\rho_+
u_+^{2}-\rho_-u_-^{2}-B(\frac{1}{\rho_+^\alpha}-\frac{1}{\rho_-^\alpha})=0.
\end{equation}

From $(\ref{5.20})$, noticing
$\widehat{\sigma^{B}}\in\Big(u_++\sqrt{ \alpha
B}\rho_+^{-\frac{\alpha+1}{2}}, u_--\sqrt{\alpha
B}\rho_-^{-\frac{\alpha+1}{2}}\Big)$, we immediately get
$(\ref{5.16})$ and $(\ref{5.17})$. The proof is finished.

\begin{Remark}
The above Lemmas \ref{lem:5.2}-\ref{lem:5.5} shows that, as
$A\rightarrow0$, the intermediate density $\rho_*^{A}$ becomes
unbounded, the velocities $\sigma_1^{A}$ and $\sigma_2^{A}$ of
shocks $S_{1}$ and $S_{2}$ and the intermediate velocity $u_{*}^{A}$
for the Riemann solutions of $(\ref{1.1})$-$(\ref{1.2})$ approach to
$\sigma^{B}$, and the intermediate density becomes a singular
measure simultaneously, which are consistent with the velocity and
the density of the $\delta$-shock solution to the generalized
Chaplygin gas equations $(\ref{1.5})$ with the same Riemann data
$(\rho_\pm,u_\pm)$ in Section 2. Thus similar to Theorem
\ref{thm:4.1}, we draw the conclusion as follows.
\end{Remark}

\begin{Theorem}\label{thm:5.1}
Let $(\rho_+,u_+)\in\rm{V}(\rho_-,u_-)$. For any fixed $A>0$, assume
that $(\rho^{A},u^{A})$ is the two-shock Riemann solution of
$(\ref{1.1})$-$(\ref{1.2})$ with Riemann data $(\rho_\pm,u_\pm)$ for
$0<A<A_{0}$ constructed in section 3. Then as $A\rightarrow0$,
$\rho^{A}$ and $\rho^{A}u^{A}$ converge in the sense of
distributions, and the limit functions of $\rho^{A}$ and
$\rho^{A}u^{A}$ are the sums of a step function and a
$\delta$-measure with weights
$$\frac{t}{\sqrt{1+(\sigma^{B})^2}}(\sigma^{B}[\rho]-[\rho u])\ \ and\ \ \frac{t}{\sqrt{1+(\sigma^{B})^2}}(\sigma^{B}[\rho u]-[\rho u^2-\frac{B}{\rho^{\alpha}}]),$$
respectively, which form a $\delta$-shock solution of $(\ref{1.5})$
with the same Riemann data $(\rho_\pm,u_\pm)$.
\end{Theorem}

\subsection{Formation of two-rarefaction-wave solutions}

Now we consider the formation of the two-rarefaction-wave
(two-contact-discontinuity) solution of $(\ref{1.1})$-$(\ref{1.2})$
with $(\ref{1.6})$ in the case $(\rho_+,u_+)\in \rm{I} (\rho_-,u_-)$
for $0<\alpha<1 (\alpha=1)$ as $A\rightarrow0$.

\begin{Lemma}\label{lem:6.7}
When $(\rho_+,u_+)\in \rm{I} (\rho_-,u_-)$, there exists a positive
parameter $A_1$ such that $(\rho_+,u_+)\in R_{1}R_{2} (\rho_-,u_-)$
when $0<A<A_{1}$.
\end{Lemma}

\noindent\textbf{Proof.} All the states $(\rho,u)$ connected with
$(\rho_-,u_-)$ by a backward shock wave $R_{1}$ or a forward shock
wave $R_{2}$ satisfy
\begin{equation}\label{5.21}
u-u_-=-\int_{\rho_{-}}^{\rho}\frac{\sqrt{An\rho^{n-1}+\frac{\alpha
B}{\rho^{\alpha+1}}}}{\rho}d\rho,\ \ \rho<\rho_{-}.
\end{equation}
or
\begin{equation}\label{5.22}
u-u_-=\int_{\rho_-}^{\rho_{}}\frac{\sqrt{An\rho^{n-1}+\frac{\alpha
B}{\rho^{\alpha+1}}}}{\rho}d\rho,\ \ \rho>\rho_{-}.
\end{equation}

When $\rho_{+}=\rho_{-}$, the conclusion is obviously true.  When
$\rho_{+}\neq\rho_{-}$, if $\rho_{+}>\rho_{-}$, by taking
$\rho>\rho_{+}$ in  $(\ref{5.22})$, we have
\begin{eqnarray}
u_+-u_-&=&\int_{\rho_-}^{\rho_{+}}\frac{\sqrt{An\rho^{n-1}+\frac{\alpha
B}{\rho^{\alpha+1}}}}{\rho}d\rho\nonumber\\
&>&\int_{\rho_-}^{\rho_{+}}\frac{\sqrt{An\rho^{n-1}}}{\rho}d\rho\nonumber\\
&=&\frac{2\sqrt{An}}{n-1}(\rho_{+}^{\frac{n-1}{2}}-\rho_{-}^{\frac{n-1}{2}})\nonumber,
\end{eqnarray}
from which we can get
\begin{equation}\label{5.23}
A<\DF{(n-1)^{2}(u_+-u_-)^{2}}{4n(\rho_{+}^{\frac{n-1}{2}}-\rho_{-}^{\frac{n-1}{2}})^{2}}.
\end{equation}
Similarly, for $\rho_{+}<\rho_{-}$, we can get the same inequality
as $(\ref{5.23})$. So we take
\begin{equation}\label{5.24}
A_{1}=\DF{(n-1)^{2}(u_+-u_-)^{2}}{4n(\rho_{+}^{\frac{n-1}{2}}-\rho_{-}^{\frac{n-1}{2}})^{2}}.
\end{equation}
The proof is finished.

When $0<A<A_{1}$, the Riemann solution of
$(\ref{1.1})$-$(\ref{1.2})$ with $(\ref{1.6})$ includes a backward
rarefaction wave $R_{1}$ and a forward rarefaction wave $R_{2}$ with
the intermediate state $(\rho_{*}^{A},u_{*}^{A})$ besides two
constant states $(\rho_{\pm},u_{\pm})$. We then have

\begin{equation}\label{5.25}R_1:\ \
\left\{\begin{array}{ll} \xi=\lambda_1^{AB}=u-\sqrt{An\rho^{n-1}+\frac{\alpha B}{\rho^{\alpha+1}}},\\
u-u_-=-\int_{\rho_{-}}^{\rho}\frac{\sqrt{An\rho^{n-1}+\frac{\alpha
B}{\rho^{\alpha+1}}}}{\rho}d\rho,\ \
\rho_{*}^{A}\leq\rho\leq\rho_{-},
\end{array}\right.
\end{equation}
and
\begin{equation}\label{5.26}R_2:\ \
\left\{\begin{array}{ll} \xi=\lambda_2^{AB}=u+\sqrt{An\rho^{n-1}+\frac{\alpha B}{\rho^{\alpha+1}}},\\
u_+-u=\int_{\rho}^{\rho_{+}}\frac{\sqrt{An\rho^{n-1}+\frac{\alpha
B}{\rho^{\alpha+1}}}}{\rho}d\rho,\ \
\rho_{*}^{A}\leq\rho\leq\rho_{+}.
\end{array}\right.
\end{equation}
Here $\rho_{*}^{A}$ is determined by
\begin{equation}\label{5.27}
u_+-u_-
=\int^{\rho_{-}}_{\rho_*^{A}}\frac{\sqrt{An\rho^{n-1}+\frac{\alpha
B}{\rho^{\alpha+1}}}}{\rho}d\rho+\int_{\rho_*^{A}}^{\rho_{+}}\frac{\sqrt{An\rho^{n-1}+\frac{\alpha
B}{\rho^{\alpha+1}}}}{\rho}d\rho
\end{equation}
Furthermore, setting
$(\rho_{*},u_{*})=\lim\limits_{A\rightarrow0}(\rho_*^{A},u_*^{A})$,
we obtain

\begin{equation}\label{5.28}
\rho_{*}^{-\frac{\alpha+1}{2}}
=\frac{(\alpha+1)(u_+-u_-)}{4\sqrt{\alpha
B}}+\frac{1}{2}(\rho_{+}^{-\frac{\alpha+1}{2}}+\rho_{-}^{-\frac{\alpha+1}{2}}),\
\
u_{*}=\frac{u_+-u_-}{2}+\frac{\sqrt{\alpha
B}}{\alpha+1}(\rho_{+}^{-\frac{\alpha+1}{2}}-\rho_{-}^{-\frac{\alpha+1}{2}}),
\end{equation}
Letting $A\rightarrow0$ in $(\ref{5.25})$ and $(\ref{5.26})$, then
for $0<\alpha<1 (\alpha=1)$, $R_{1}$ and $R_{2}$ become the backward
rarefaction wave (the backward contact discontinuity)
$\overleftarrow{R }$ and the forward rarefaction wave (the forward
contact discontinuity) $\overrightarrow{R}$, respectively, as
follows:
\begin{equation}\label{5.29}
 \overleftarrow{R}:\left\{\begin{array}{ll} \xi=\lambda_1^B=u-\sqrt{\alpha B}\rho_{}^{-\frac{\alpha+1}{2}},\\
u-\frac{2\sqrt{\alpha
B}}{1+\alpha}\rho_{}^{-\frac{\alpha+1}{2}}=u_--\frac{2\sqrt{\alpha
B}}{1+\alpha}\rho_{-}^{-\frac{\alpha+1}{2}},\ \ \rho_->\rho>\rho_*,
\end{array}\right.
\end{equation}
and
\begin{equation}\label{5.30}
\overrightarrow{R}:\left\{\begin{array}{ll} \xi=\lambda_2^B=u+\sqrt{\alpha B}\rho_{}^{-\frac{\alpha+1}{2}},\\
u+\frac{2\sqrt{\alpha
B}}{1+\alpha}\rho_{}^{-\frac{\alpha+1}{2}}=u_++\frac{2\sqrt{\alpha
B}}{1+\alpha}\rho_{+}^{-\frac{\alpha+1}{2}},\ \ \rho_+>\rho>\rho_*.
\end{array}\right.
\end{equation}

As a conclusion, for the case $(\rho_+,u_+)\in \rm{I} (\rho_-,u_-)$
,as $A\rightarrow0$, the two rarefaction wave $R_{1}$ and $R_{2}$ in
$(\ref{5.25})$ and $(\ref{5.26})$ approach the two rarefaction waves
(contact discontinuities) $\overleftarrow{R }$ and
$(\overrightarrow{R})$ in $(\ref{5.29})$ and $(\ref{5.30})$ for
$0<\alpha<1 (\alpha=1)$, and the intermediate state
$(\rho_*^{A},u_*^{A})$ tends to the state $(\rho_*,u_*)$ in
$(\ref{5.28})$. In summary, in this case, we have the following
result.

\begin{Theorem}\label{thm:5.2}
Let $(\rho_+,u_+)\in\rm{I}(\rho_-,u_-)$. For any fixed $A>0$, assume
that $(\rho^{A},u^{A})$ is the two-shock Riemann solution of
$(\ref{1.1})$-$(\ref{1.2})$ with Riemann data $(\rho_\pm,u_\pm)$ for
$0<A<A_{1}$ constructed in section 3. Then as $A\rightarrow0$, the
limit of the Riemann solution
 $(\rho^{A}, u^{A})$ is two rarefaction waves (contact
discontinuities) connecting the constant states $(\rho_\pm,u_\pm)$
and the intermediate nonvacuum state as follows:
\begin{equation}\nonumber
\lim\limits_{A\rightarrow0}(\rho,u)(\xi)=\left\{\begin{array}{ll}
(\rho_-,u_-),\ \ \ \ \ \ -\infty<\xi\leq u_{-}-\sqrt{\alpha B}\rho_{-}^{-\frac{\alpha+1}{2}},\\
(\rho_*,u_*),\ \ \ \ \ \ \ \ u_--\sqrt{\alpha B}\rho_{-}^{-\frac{\alpha+1}{2}}\leq\xi\leq u_++\sqrt{\alpha B}\rho_{+}^{-\frac{\alpha+1}{2}},\\
(\rho_+,u_+),\ \ \ \ \ \ \ u_++\sqrt{\alpha
B}\rho_{+}^{-\frac{\alpha+1}{2}}\leq\xi<\infty,
\end{array}\right.
\end{equation}
which is exactly the Riemann solution to the (generalized) Chaplygin
gas equations $(\ref{1.5})$ with the same Riemann data
$(\rho_\pm,u_\pm)$ for $0<\alpha<1 (\alpha=1)$.
\end{Theorem}

\section{Conclusions and discussions}
In this paper, we have consider two kinds of the flux approximation
limit of Riemann solutions to extended Chaplygin gas equations and
studied the concentration and the formation of delta shock during
the limit process. Moreover, we have proved that the vanishing
pressure limit of the Riemann solutions to extended Chaplygin gas
equations is just the corresponding ones to trasport equations, and
when extended Chaplygin pressure approaches the generalized
Chaplygin pressure, the limit of the Riemann solutions to extended
Chaplygin gas equations is just the corresponding ones to the
generalized Chaplygin gas equations. In fact, one can further prove
that when the extended Chaplygin pressure approaches the pressure
for the perfect fluid, i.e., $B\rightarrow0$ for fixed $A$, the
limit of the Riemann solutions to the extended Chaplygin gas
equations is just the corresponding ones to the Euler equations for
perfect fluids.

On the other hand, recently, Shen and Sun have studied the Riemann
problem for the nonhomogeneous tranport equations, and the
nonhomogeneous (generalized) Chaplygin gas equations with
coulomb-like friction, see \cite{Shen1,Shen2,Sun}. Similarly, we
will also consider the Riemann problem for the nonhomogeneous
extended Chaplygin gas equations with coulomb-like friction.
Furthermore, we will consider its flux approximation limit and
analyze the relations Riemann solutions among the nonhomogeneous
extended Chaplygin gas equations, the generalized Chaplygin gas
equations and the nonhomogeneous trasport equations. These will be
left for our future work.

 \bigbreak


\bigbreak

\end{document}